\documentclass{article}
\usepackage{geometry} 
\usepackage{amsmath, amsthm, amssymb, mathtools}
\usepackage{color}
\newtheorem{lemma}{Lemma}[section]

\newtheorem{proposition}[lemma]{Proposition}

\newtheorem{theorem}[lemma]{Theorem}

\theoremstyle{definition}
\newtheorem{remark}[lemma]{Remark}
\newcommand{\e}{\varepsilon} 

\newcommand{\R}{\mathbb{R}}
\newcommand{\weak}{\rightharpoonup}

\newcommand{\mres}{\mathbin{\vrule height 1.6ex depth 0pt width
0.13ex\vrule height 0.13ex depth 0pt width 1.3ex}}

\begin{document} 
\title{Higher-order singular perturbation models\\
 for phase transitions}
\author{ Giuseppe Cosma Brusca, Davide Donati, \\
{\small SISSA}\\ \small
via Bonomea 265\\ \small
34146 Trieste, Italy\\
\and
Margherita Solci\footnote{email: {\tt margherita@uniss.it}}\\ 
{\small DADU, Universit\`a di Sassari}\\ {\small piazza Duomo 6}\\{\small 07041 Alghero, Italy}}
\date{}							
\maketitle
\begin{abstract} Variational models of phase transitions take into account double-well energies singularly perturbed by gradient terms, such as the Cahn-Hilliard free energy. The derivation by $\Gamma$-convergence of a sharp-interface limit for such energy is a classical result by Modica and Mortola. We consider a singular perturbation of a double-well energy by derivatives of order $k$, and show that we still can describe the limit as in the case $k=1$ with a suitable interfacial energy density, in accord with the case $k=1$ and with the case $k=2$  previously analyzed by Fonseca and Mantegazza. The main isssue is the derivation of an optimal-profile problem on the real line describing the interfacial energy density, which must be conveniently approximated by minimum problems on finite intervals with homogeneous condition on the derivatives at the endpoints up to order $k-1$. To that end a careful study must be carried on of sets where sequences of functions with equibounded energy are ``close to the wells'' and have ``small derivatives'', in terms of interpolation inequalities and energy estimates.

\smallskip



\end{abstract}
\section{Introduction}
A widely used model for phase transitions can be formulated as a minimum principle for double-well energies perturbed by a gradient term, the so-called Cahn-Hilliard energy
\begin{equation}\label{1}
\int_\Omega W(u)\,dx+\e^2\int_\Omega |\nabla u|^2\,dx,
\end{equation}
where $\Omega$ is a bounded open set in $\mathbb R^d$ with Lipschitz boundary, $u\colon\Omega\to\mathbb R$, and
$W$ is a non-negative function with only two zeroes; say, $-1$ and $1$.
The variable $u$ is usually subjected to an additional integral constraint, which forbids trivial minimizers.
The asymptotic analysis of such energies can be stated as the computation of a $\Gamma$-limit.
To that end, the energies have to be scaled by $\frac{1}{\e}$ in order not to have a trivial result.
In a seminal paper \cite{MM}, Modica and Mortola have computed this $\Gamma$-limit and shown that its domain is
the set $BV(\Omega;\{-1,1\})$ of the functions of bounded variations $u$ which only take the values $1$ and $-1$.  
In terms of the set of essential discontinuity points $S(u)$, which can be interpreted as the interface between the 
regions where $u=1$ and $u=-1$, the $\Gamma$-limit takes the form
\begin{equation}\label{2}
m_1\mathcal H^{d-1}(S(u)\cap\Omega)=m_1{\rm Per}(\{u=1\},\Omega),
\end{equation}
where $m_1$ is the value given by the {\em optimal-profile problem} 
\begin{equation}\nonumber
m_1=\inf \Bigl\{\int_{-\infty}^{+\infty} \big(W(u)+(u')^2\big)\, dt: u\in H^1_{\textup{loc}}(\mathbb{R}), \lim_{t\to\pm\infty} u(t)=\pm1\Bigr\},
\end{equation}
$\mathcal H^{d-1}$ is the $(d-1)$-dimensional Hausdorff measure and Per$(E,\Omega)$ 
denotes the perimeter of a set $E$ in $\Omega$ in the sense of sets of finite perimeter \cite{Maggi}. 
Hence, after showing that the integral constraint for $u$ is compatible with $\Gamma$-convergence and rereading it as a constraint for the measure $|\{u=1\}|$, the Fundamental Theorem of $\Gamma$-convergence allows one to interpret the limit of minimizers as a set of minimal perimeter, thus confirming a {\em minimal-interface criterion} as conjectured by Gurtin \cite{Gurtin}.  
It is worth remarking that the problem defining $m_1$ is set on the real line, and depends only on $W$. This is due to the fact that minimizing functions at given $\e$ have a one-dimensional structure, following scaled optimal profiles in the direction orthogonal to the limit interface. This highlights that in fact the relevant analysis is in dimension $d=1$, the higher-dimensional case being derived from that one using by-now-standard technical tools \cite{BLN98}.  

In this paper, we describe the behaviour of the functionals in \eqref{1} when the gradient is substituted by a suitably-scaled higher-order derivative; more precisely, of functionals
\begin{equation}\label{4}
\int_\Omega W(u)\,dx+\e^{2k}\int_\Omega \|D^{k} u\|^2\,dx.
\end{equation}
Besides the motivations related to the fine study of models of phase transition, for which we refer to the detailed introduction of the paper by Fonseca and Mantegazza \cite{Fonseca2000}, where they analyze the case $k=2$, higher-order perturbation of the form \eqref{4} have also been used to approximate free-discontinuity energies (see \cite{ABG,BDS,mo}), to the study of the Perona-Malik model in Visual Reconstruction (see \cite{bcg,gp}), to phase-field models for free-discontinuity problems (see e.g.~\cite{Bach,FLL}), to improve the regularity of the associated phase-field variable (see e.g.~\cite{BBL}) and provide numerical approximation for high-order phase-field models for Brittle Fracture \cite{Negri}.
In a different but related context, higher-order perturbations have also been used in the study of multi-well problems in elasticity (see e.g.~\cite{CFL,Cos,Ho,ChCo}) and for higher-order materials (see e.g~\cite{CDMFL}), and have been examined in relation to other limit procedures, such as dimension reduction (see e.g.~\cite{IZ}).

Here, we consider functionals \eqref{4} for arbitrary $k\ge 3$, for which the relevant analysis by $\Gamma$-convergence is that of the scaled one-dimensional functionals
\begin{equation}\label{5}
\frac1\e\int_I W(u)\,dx+\e^{2k-1}\int_I (u^{(k)})^2\,dx,
\end{equation}
defined for $u\in H^k(I)$, where $I$ is an open interval. We show that their $\Gamma$-limit in the $L^1(I)$ topology and with respect to the convergence in measure is still of the form 
\begin{equation}\nonumber
m_k\#(S(u)\cap I),\qquad u\in BV(I;\{-1,1\}),
\end{equation}
where $m_k$ is the value given by the optimal-profile problem 
\begin{equation}\nonumber
m_k=\inf \Bigl\{\int_{-\infty}^{+\infty} \big(W(u)+(u^{(k)})^2\big)\, dt: u\in H^k_{\textup{loc}}(\mathbb{R}), \lim_{t\to\pm\infty} u(t)=\pm1\Bigr\},
\end{equation}
in accord with the results in \cite{MM}  for $k=1$ and \cite{Fonseca2000} for $k=2$. 

The case $k>2$ requires a delicate analysis of  sequences $u_\e$ with equibouded energy.
To that end, as for $k\le 2$, it is necessary to analyze the property of the sets
\begin{eqnarray}\nonumber
\{ t\in I:  |u_\e(t)-1|<\eta\},\qquad
\{ t\in I:  |u_\e(t)+1|<\eta\};
\end{eqnarray}
that is, of the set of points with $u_\e$ close to the bottom of the wells.
In the case $k=1$ these sets are described from properties of the only part with the potential $W$ in the energy. The energetic behaviour of the sequence can then be described through the approximate optimal-profile problems
\begin{eqnarray}\nonumber
m_1(\eta,T)=\inf \Bigl\{\int_{-T}^{T} \big(W(u)+(u^{(k)})^2\big)\, dt: u\in H^k(-T,T),\\\nonumber
 |u(-T)+1|<\eta,
|u(T)-1|<\eta\Bigr\}.
\end{eqnarray}
Heuristically, once it is remarked that for $\eta$ sufficiently small $m_1(\eta,T)\geq c>0$ for all $T$, every transition of $u_\e$ between values close to $1$ to values close to $-1$ or the converse requires a fixed amount of energy, from which we can conclude the pre-compactness of the sequence. By optimizing in $T$ and $\eta$ we can finally prove the lower estimate with $m_1$.

In the case $k=2$ Fonseca and Mantegazza \cite{Fonseca2000} have noticed that it is important to 
study points with $u_\e$ close to the bottom of the wells and of small derivative.
They have studied the set of such points using a bound on the $L^{4/3}$-norm of $u'_\e$ obtained by an interpolation inequality between the $L^1$-norm of $u_\e$ and the $L^2$-norm of $u''_\e$. The properties of such points allow one to carry on the argument above with 
\begin{eqnarray}\nonumber
m_2(\eta,T)=\inf \Bigl\{\int_{-T}^{T} \big(W(u)+(u^{(k)})^2\big)\, dt: u\in H^k(-T,T),\\ \nonumber |u(-T)+1|<\eta, 
|u(T)-1|<\eta, |u^\prime(\pm T)|<\eta\Bigr\}
\end{eqnarray}
in the place of $m_1(\eta,T)$.
As it is remarked in \cite{Fonseca2000}, that interpolation argument requires at least linear growth for $W$ at infinity, and is not easily extended to $k\ge 3$.

In the present work we use a different and direct strategy, which in particular does not use any growth hypothesis on $W$. Its main point is the analysis of the properties of the sets
\begin{eqnarray} \nonumber
A_\e^{\eta,+}=\{ t\in I:  |u_\e(t)-1|<\eta, |u_\e^{(\ell)}(t)|<\e^\ell\eta\hbox{ for all }\ell\le k-1\},
\end{eqnarray}
and the corresponding $A_\e^{\eta,-}$ for which the first condition is replaced by $|u_\e(t)+1|<\eta$; that is,
of the set of points with $u_\e$ close to the bottom of the wells and with small derivatives up to  order $k-1$.
The main feature of such sets is that all elements of any family of intervals $I_\e$ of length much larger than $\e$ contain points of either 
$A_\e^{\eta,+}$ or $A_\e^{\eta,-}$. This is obtained by using $L^2$ estimates on all derivatives $u^{(\ell)}_\e$ for $\ell\in\{1,\ldots,k-1\}$, which are achieved by simultaneously using interpolations between the $L^2$-norm of $u_\e-1$ (or $u_\e+1$) and the $L^2$-norm of $u^{(k)}_\e$. The argument is an adaptation of a technique recently devised by Solci \cite{Solci2024} and used to obtain free-discontinuity functionals by singular perturbation,
generalizing to higher-order derivatives a result by Alicandro et al.~\cite{ABG} (see also Bouchitt\'e et al \cite{BDS}). In the present situation the argument is simpler since in the limit we may only have discontinuities with a jump size equal to $2$, and not of arbitrary size. 
Note that the use of the $L^2$-norm of $u_\e\pm1$ does not require any growth condition on $W$ since it is computed only on intervals in $A_\e^{\eta,\pm}$, on which we have an $L^\infty$ bound on $u_\e\pm1$.
 The argument for $k=1$ and $k=2$ is then carried over by the use of the approximate optimal-profile problems
\begin{eqnarray}\nonumber
m_k(\eta,T)=\inf \Bigl\{\int_{-T}^{T} \big(W(u)+(u^{(k)})^2\big)\, dt: u\in H^k(-T,T), |u(-T)+1|<\eta, \\\nonumber
|u(T)-1|<\eta, |u^{(\ell)}(\pm T)|<\eta \hbox{ for all }\ell\le k-1\Bigr\},
\end{eqnarray}
which lead to $m_k$. An additional issue with respect to the cases $k=1$ and $k=2$ here is to show that transitions between the bottom of the wells where we do not have a bound on the derivatives at the endpoints are negligible.

We finally note that the analysis of functionals \eqref{4} can be reduced to the one-dimensional case by applying a Fubini argument and considering one-dimensional sections, once we take $\|D^ku\|$ as the operator norm on the tensor $D^ku$.
As in the cases $k=1$ and $k=2$ the $\Gamma$-limit has the form \eqref{2} with $m_1$ replaced by $m_k$.

\section{Statement of the result}
Let $W\colon\mathbb R\to [0,+\infty)$ be a continuous function such that
$W(-1)=W(1)=0$ and  
$$W(z)\geq \alpha_W \min\{(z+1)^2, (z-1)^2, \beta_W\}$$ 
for some $\alpha_W, \beta_W>0$. 

Let $k\geq 1$ be an integer and let $(a,b)$ be an open bounded interval. 
For any $\e>0$ we define on $H^k(a,b)$ the functional $F_\e$ given by  
\begin{equation}\label{defFe}
F_\e(u)=\int_{(a,b)}\Big(\frac{1}{\e}W(u)+ \e^{2k-1}(u^{(k)})^2\Big)\, dt
\end{equation} 
with $u\in H^k(a,b)$. 

\smallskip 

We will proof the following compactness and $\Gamma$-convergence theorem,  noting that the result if $k=1$ and $k=2$ has been shown in \cite{MM} and in \cite{Fonseca2000}, respectively. 
{
Note that the domain of the $\Gamma$-limit in this case is the space $BV((a,b);\{-1,1\})$ of functions with bounded variation $u$ on the interval $(a,b)$ with $u(t)\in \{-1,1\}$ almost everywhere, which are just the piecewise-constant functions on $(a,b)$ with values equal to $-1$ or $1$. 
}

\begin{theorem}[Coerciveness and $\Gamma$-convergence]\label{main} 
Let $F_\e$ be defined by \eqref{defFe} on $H^k(a,b)$. Then, the family $\{F_\e\}$ is equicoercive with respect to the convergence in measure; that is, 
for all $\e_j\to 0$ and $\{u_j\}$ with $\sup_j F_{\e_j}(u_j)<+\infty$, there exists $u\in BV((a,b);\{-1,1\})$ such that, up to subsequences, $u_j\to u$ in measure. 
Moreover, the family $\{F_\e\}$ $\Gamma$-converges with respect to the convergence in measure and in $L^1(a,b)$ to the functional $F$ with domain $BV((a,b);\{-1,1\})$ given by 
$$F(u)=m_k \# S(u),$$
where 
\begin{equation}\label{defmklim}
    m_k = \inf \Bigl\{\int_{-\infty}^{+\infty} \bigl(W(u)+(u^{(k)})^2\bigr)\, dt: u\in H^k_{\textup{loc}}(\mathbb{R}), \lim_{t\to\pm\infty} u(t)=\pm1\Bigr\}
\end{equation}
and $S(u)$ denotes the jump set of $u$. 
More precisely, for all $u\in BV((a,b);\{-1,1\})$ we have that 

\smallskip 

{\rm (i)} for all sequences $u_\e\to u$ in measure,  
\begin{equation}\label{liminf1}\nonumber
\liminf_{\e\to 0}F_\e(u_\e)\geq m_k \#S(u); 
\end{equation} 

\smallskip 

{\rm (ii)} there exists a recovery sequence $u_\e\to u$ in $L^1(a,b)$ such that   
\begin{equation}\label{liminf2}\nonumber
\limsup_{\e\to 0}F_\e(u_\e)=m_k \#S(u). 
\end{equation} 
\end{theorem}

{

\bigskip

From the one-dimensional result, we will recover an analog in dimension $d>1$ applying blow-up and localization techniques. We consider $\Omega\subset\mathbb{R}^d$ an open and bounded set with Lipschitz boundary, and, for $\e>0$, the functional defined for $u\in H^k(\Omega)$ as
\begin{equation}\label{eq:def F multidim}
        F_\e(u):=\int_\Omega\Big(\frac{1}{\e}W(u)+\e^{2k-1}\|D^ku\|^2\Big) dx,
\end{equation}
where $W\colon\R\to [0,+\infty)$ is a function as above and $\|D^ku\|$ denotes the operator norm of the tensor $D^ku$. In accord with Theorem \ref{main}, the domain of the $\Gamma$-limit will be the space $BV(\Omega;\{-1,1\})$ of functions $u$ with bounded variation on the set $\Omega$ with $u(x)\in \{-1,1\}$ almost everywhere. We recall that in this case the set of essential discontinuity points $S(u)$ is a set with finite $(d-1)$-dimensional Hausdorff measure, and that  $\mathcal H^{d-1}(S(u))$ is equal to the perimeter of the set $\{x: u(x)=1\}$ in $\Omega$ in the sense of the theory of sets of finite perimeter (see e.g.~\cite{Maggi,BLN98}).

\begin{theorem}\label{main-d}
The family $\{F_\e\}$ defined by \eqref{eq:def F multidim} on $H^k(\Omega)$ $\Gamma$-converges with respect to the  $L^1(\Omega)$-convergence to the functional $F$ with domain $BV((\Omega);\{-1,1\})$ given by 
$$F(u)=m_k \mathcal{H}^{d-1}(S(u))=m_k{\rm Per}(\{u=1\};\Omega),$$ where $m_k$ is as defined in \eqref{defmklim}, $S(u)$ denotes the jump set of $u$, and {\rm Per}$(A;\Omega)$ the perimeter of $A$ in $\Omega$ in the sense of the theory of sets of finite perimeter.     
\end{theorem}

\begin{remark} Note that, if we strengthen the growth assumptions on the double-well potential $W$ or alternatively we assume the equi-integrability of $\{u_\e\}$, we may deduce that a sequence $\{u_\e\}$ such that $\sup_\e F_\e(u_\e)<+\infty$ is strongly precompact in $L^1(\Omega)$. This follows by a result due to Alberti, Bouchitt\`e, Seppecher \cite[Theorem 6.6]{alberti1998phase}.
\end{remark}
}
To prove the results, it is useful to introduce the localized functionals 
\begin{equation}\label{defFeloc}\nonumber
F_\e(u;I)=\int_{I}\Big(\frac{1}{\e}W(u)+ \e^{2k-1}(u^{(k)})^2\Big)\, dt
\end{equation}
for any interval $I\subset (a,b)$. Finally, to simplify the notation, in the following we only consider $(a,b)=(0,1)$, the changes in the proofs for a general interval being obvious.

\section{Analysis of oscillations close to $1$ and $-1$} 
In the proofs of the compactness and of the lower bound, we will make use of 
some optimal-profile problems in order to estimate the number of transitions between wells of sequence of functions with equibounded energy. 
In particular, we use the family of minimum problems given by 
\begin{eqnarray}\label{defment}
   && \hskip-5mm  \nonumber m_k(\eta, N,T):=\inf \Bigl\{\int_{-T}^T \bigl(W(u)+(u^{(k)})^2\bigr)\, dt: u\in H^k(-T,T),\nonumber \\ 
   && \hskip3cm  |u(-T)+1|\leq \eta, 
    |u(T)-1|\leq \eta,\nonumber \\
    && \hskip3cm  |u^{(\ell)}(\pm T)|\leq\textstyle\frac{1}{N} \,\,\textup{for all } \ell\in\{1,..,k-1\}\Bigr\},
\end{eqnarray}
describing transitions between points where the functions are close to  $-1$ or $1$ and their derivatives up to order $k-1$ are ``small''.

In this section, we will study the properties of intervals of points where functions $u_\e$ with equibounded energy are close to either $-1$ or $1$, and their derivatives up to order $k-1$ are small. 
We will show that every sufficiently large interval (with respect 
to $\e$) intersects such sets. The technical tool used to prove this result is the simultaneous application of interpolation inequalities for all derivatives in the intervals where $u_\e$ is close to either $-1$ or $1$. This is an adaptation and simplification of the technique devised in \cite{Solci2024} to treat an approximation of free-discontinuity problems by singular perturbation with higher-order derivatives.

\smallskip 

Let $\{u_\e\}$ be a family of $H^k(0,1)$ functions with equibounded energy, and let 
\begin{equation}\label{defS}\nonumber
S=\sup_{\e>0} F_\e(u_\e). 
\end{equation}
For any fixed $\eta\in(0,1)$ and $\e>0$, we define the set 
\begin{equation}\label{defae}
A_\e^\eta=\{t\in(0,1): ||u_\e|-1|<\eta\}.
\end{equation} 
Recalling that 
$W(z)\geq\alpha_W\min\{(z-1)^2,(z+1)^2,\beta_W\}$ for some $\alpha_W,\beta_W>0$, 
we will consider $\eta\in(0,\sqrt\beta_W)$, so that $W(z)\geq \alpha_W(z-1)^2$ if $|z-1|\leq \eta$ and $W(z)\geq \alpha_W(z+1)^2$ if $|z+1|\leq \eta$. 

We now fix an integer $N\geq 1$, and prove an upper estimate for the length of an interval $I$ such that at all points of $I\cap A_\e^\eta$ there exists at least one $\ell\in\{1,\dots, k-1\}$ such that $|u_\e^{(\ell)}|\geq \frac{1}{N\e^{\ell}}$.  
\begin{lemma}\label{stimaint} 
Let $\eta\in (0,\sqrt{\beta_W})$ and $N\geq 1$ be fixed, and let $\{u_\e\}$ be such that $\sup_\e F_\e(u_\e)=S<+\infty$. 
Then, there exists a constant $R_k(\eta, N, S)>0$ such that, for all $\e>0$ and for any interval $I$ satisfying the hypotheses: 
\begin{enumerate}
\item[{\rm (i)}] $I\subset A_\e^\eta$; 
\item[{\rm (ii)}] $\displaystyle\Big|\Big\{t\in I: |u_\e^{(\ell)}(t)|< \frac{1}{N\e^{\ell}} \ \hbox{\rm for all } \ell\in\{1,\dots, k-1\}\Big\}\Big|=0$;  
\end{enumerate} 
it follows that 
$$|I|\leq R_k(\eta, N, S) \e.$$ 
\end{lemma} 
\begin{proof} 
The proof is based on a classical interpolation estimate. 
There exists a constant $r_k>0$ such that for any interval $I$ and for any $v\in H^k(I)$ 
the following estimate holds
\begin{equation}\label{interpolation1} 
\|v^{(\ell)}\|_{L^2(I)}\leq r_k \Big(\|v\|^\theta_{L_2(I)} \|v^{(k)}\|^{1-\theta}_{L_2(I)} + |I|^{-\ell} \|v\|_{L_2(I)}\Big) 
\end{equation}
for all $\ell\in\{1,\dots, k-1\}$, with $\theta=\frac{k-\ell}{k}$ (see e.~g. \cite[Theorem 7.41]{leoni}).  
Setting 
$$a=\e^{-\theta}\|v\|_{L^2(I)}^{2\theta} \ \ \hbox{\rm and } \ \ b=\e^{(2k-1)(1-\theta)}\|v^{(k)}\|_{L^2(I)}^{2(1-\theta)},$$ 
by using the convexity inequality  
$ab\leq \frac{1}{p} a^p+\frac{1}{q} b^q$   
with $p=\frac{1}{\theta}$ and $q=\frac{1}{1-\theta}$, 
we then deduce that there exists a constant $R_k>0$ such that 
\begin{equation}\label{interpolation2} 
\e^{2\ell-1}\|v^{(\ell)}\|^2_{L^2(I)}\leq R_k \Big(\frac{1}{\e}\|v\|^2_{L_2(I)} +\e^{2k-1}\|v^{(k)}\|^2_{L_2(I)} + \frac{\e^{2\ell-1}}{|I|^{2\ell}} \|v\|^2_{L_2(I)}\Big). 
\end{equation}

Let now $I$ be as in the hypothesis. It is not restrictive to suppose that $|u_\e-1|\leq \eta$ in $I$, the case $|u_\e+1|\leq \eta$ being completely analogous. Choosing $v=u_\e-1$ in \eqref{interpolation2}, and recalling that $W(u_\e)\geq \alpha_W(u_\e-1)^2$ in the interval $I$, we get 
\begin{equation}
\label{interpolation3} 
\e^{2\ell-1}\|u_\e^{(\ell)}\|^2_{L^2(I)}\leq \widetilde R_k \Big(F_\e(u_\e; I) + \e^{2\ell-1}|I|^{-2\ell} \|u_\e-1\|^2_{L_2(I)}\Big),  
\end{equation} 
for any $\e>0$, where $\widetilde R_k=R_k \max\{\frac{1}{\alpha_W},1\}$. 

Note that hypothesis (ii) implies that for at least one $\ell\in\{1,\dots, k-1\}$ 
$$\Big|\Big\{t\in I: |u_\e^{(\ell)}(t)|\geq \frac{1}{N\e^{\ell}}\Big\}\Big|\geq \frac{|I|}{k+1}.$$ 
By using this lower estimate in \eqref{interpolation3}, recalling that $(u_\e-1)^2\leq \eta^2$ and that $F_\e(u_\e; I)\leq S$ we obtain 
\begin{equation}
\label{interpolation4} \nonumber
\e^{-1}\frac{1}{(k+1) N^2} |I| \leq \widetilde R_k \Big(S + \eta^2 \e^{2\ell-1}|I|^{1-2\ell}\Big)
\end{equation} 
for any $\e>0$.
We then conclude that 
 \begin{equation}
 \label{interpolation5} \nonumber
|I| \leq \widetilde R_k (k+1) N^2 S \e + \widetilde R_k (k+1) N^2 \eta^2 \e^{2\ell}|I|^{1-2\ell}. 
\end{equation} 
Hence, if $|I| > 2 \widetilde R_k (k+1) N^2 S \e$, it follows that 
$$|I| \leq  2\widetilde R_k (k+1) N^2 \eta^2 \e^{2\ell}|I|^{1-2\ell},$$  
so that $|I|^{2\ell} \leq  2 \widetilde R_k (k+1) N^2 \eta^2 \e^{2\ell}$. Hence, 
$$|I|\leq \max\Big\{2 \widetilde R_k (k+1) N^2 S, \max_{1\leq\ell\leq k-1}\big\{(2\widetilde R_k (k+1) N^2 \eta^2)^{\frac{1}{2\ell}}\big\}\Big\}\e.$$ 
With a little abuse of notation, by letting $R_k(\eta, N, S)$ denote the constant in this inequality, we obtain the claim. 
\end{proof}

We set 
\begin{equation}\label{defaeN} 
A_\e^\eta(N)=\Big\{t\in A_\e^\eta: |u_\e^{(\ell)}(t)|< \frac{1}{N\e^{\ell}} \ \hbox{\rm for all } \ell\in\{1,\dots, k-1\}\Big\},
\end{equation}
where we recall that $A_\e^\eta$ is defined in \eqref{defae} as the set of $t$ such that $||u_\e(t)|-1|\leq \eta$.  
The sets $A_\e^\eta(N)$ describe the regions close to $\{\pm1\}$ where all derivatives of $u_\e$ up to $k-1$ are ``small''. We now show that points in such sets can not have distance much larger than $\e$. 

\begin{lemma}\label{lunghint} 
Let $\eta\in (0,\sqrt{\beta_W})$, $N\geq 1$ and let $\{u_\e\}$ be such that $\sup_\e F_\e(u_\e)=S<+\infty$. 
Let $\e_r\to 0$ as $r\to+\infty$, and 
let $\{I_{\e_r}\}$ be a family of intervals such that 
$$\lim_{r\to+\infty}\frac{|I_{\e_r}|}{\e_r}=+\infty.$$
Then, there exists $r_0\in \mathbb N$ such that 
$$|I_{\e_r}\cap A^\eta_{\e_r}(N)|>0$$
for all $r\geq r_0$. 
\end{lemma}

\begin{proof} 
Let $\{I_{\e_r}\}$ be the family of intervals in the hypothesis. 
By contradiction, we suppose that 
$|I_{\e_r}\cap A^\eta_{\e_r}(N)|=0$ 
for a subsequence (not relabeled) $\e_r\to 0$ as $r\to+\infty$. 
The proof of the contradiction is divided in two steps. 

\medskip 

\noindent{\em Step 1.} Let $r$ be fixed, and let $J\subset I_{\e_r}$ be an interval such that 
$$|J|\geq \e_r \Big(\frac{2S}{\eta^2} + (2k+3) R_k(\eta, N,S)\Big),$$
where $R_k(\eta, N,S)>0$ is the constant given in Lemma \ref{stimaint}. 
Then for all $\ell\in\{1,\dots, k-1\}$ there exists a point in $J$ where $u_{\e_r}^{(\ell)}$ vanishes.  
\smallskip 

We first note that $|\{||u_{\e_r}|-1|\geq \eta\}|\leq \e_r\frac{S}{\eta^2}$. Since $|I_{\e_r}\cap A^\eta_{\e_r}|=0$, we can apply Lemma \ref{stimaint} deducing that $J$ contains at least $2k+3$ intervals 
such that in each one of them $|u_{\e_r}-1|<\eta$, with the equality holding at the endpoints,
or $|u_{\e_r}+1|<\eta$, again with the equality holding at the endpoints. 
Let $\{I_n=(a_n, b_n)\}_{n=1}^{2k+3}$ denote the set of these intervals, with $a_n$ and $b_n$ increasing. 

Now, we show that in each interval $(a_n,a_{n+2})$ with $n\leq 2k+1$ there exists at least a point where $u_{\e_r}^\prime$ vanishes. Indeed, let $n\leq 2k+1$ and suppose that $|u_{\e_r}-1|<\eta$ in $I_n$. If   $u_{\e_r}^\prime\neq 0$ in the whole $I_n$, then $u_{\e_r}$ is strictly monotone. 
If $u_{\e_r}$ is strictly increasing in $I_{n}$, there exists a local maximizer in $(b_n,a_{n+1})$. Otherwise, $u_{\e_r}$ can not be decreasing in $(b_n, a_{n+2})$, so that we it contains a local minimizer, proving the claim.  
The same argument can be applied if $|u_{\e_r}+1|<\eta$ in $I_n$. 

Hence, there exist at least $k$ internal points in $J$ where $u_{\e_r}^\prime=0$. By an iteration of the application of Lagrange Theorem, this implies that for each $\ell\in\{1,\dots, k-1\}$ there exists at least a point where $u_{\e_r}^{(\ell)}=0$. 

\medskip 

\noindent{\em Step 2.}    
Let $N_r$ be the integer part of $\frac{|I_{\e_r}|}{\e_r L}$, 
where 
$$L=L_k(\eta, N, S)=\frac{2S}{\eta^2} + (2k+3) R_k(\eta, N,S),$$  
so that each interval $I_{\e_r}$ contains $N_r$ disjoint intervals with length $\e_r L$ 
and by hypothesis $N_r\to+\infty$. 
Then, for any $r$ there exists one of these intervals, denoted by $J_r=(a_r,a_r+\e_r L)$, such that 
$$F_{\e_r}(u_{\e_r};J_r)\leq \frac{S}{N_r}.$$ 
We can then write 
\begin{eqnarray}\label{stimav} 
\frac{S}{N_r}&\geq&\int_{a_r}^{a_r+\e_r L}\Big(\frac{1}{\e_r}W(u_{\e_r})+\e_r^{2k-1}(u_{\e_r}^{(k)})^2\Big)\, dt\nonumber\\
&=&\int_{0}^{L}\Big(W(v_{r})+(v_{r}^{(k)})^2\Big)\, ds, 
\end{eqnarray} 
where $v_r(t)=u_{\e_r}(\e_r t+a_r)$. 
By Step 1, for any $\ell\in\{1,\dots, k-1\}$ there exists a point $t_r^\ell$ in $J_r$ where $u_{\e_r}^{(\ell)}$ vanishes, hence $v_r^{(\ell)}=0$ in $s_r^\ell=\frac{t_r^\ell-a_r}{\e_r}\in (0, L)$. 
An iterated application of the Fundamental Theorem of Calculus then gives 
$$\|v_r\|_{H^{k}(0,L)}\leq C(L) \big(\min_{[0,L]} |v_r|^2+\|v_r^{(k)}\|^2_{L^2(0,L)}\big),$$
where $C(L)>0$ depends on $L$ (and then on $\eta, N, S$ and $k$). Moreover, note that there exists at least a point in $J_r$ such that $||u_{\e_r}|-1|=\eta$, and the same holds for $v_r$ in a point in $(0,L)$. 
Then, since $N_r\to +\infty$, by \eqref{stimav} we deduce that $v_r$ is equibounded in $H^k(0,L)$, thus $v_r\weak v$ in $H^k(0,L)$ as $r\to+\infty$. 
Since $v_r\to v$ uniformly, by \eqref{stimav} we deduce that 
$$0=\liminf_{r\to+\infty}\int_{0}^{L}W(v_{r})\, ds\geq \int_{0}^{L}W(v)\, ds,$$ 
so that either $v= 1$ in $(0,L)$ or $v= -1$ in $(0,L)$. This gives a contradiction since 
for any $r$ there exists a point $s_r\in(0,L)$ such that 
$||v_{\e_r}(s_r)|-1|=\eta$. 
\end{proof}

\section{Optimal-profile problems}
In the proof of the equicoerciveness of the family $\{F_\e\}$ and of the lower estimate for the $\Gamma$-limit, we will make use of the optimal-profile problem 
\begin{equation}\label{defmklim-2}\nonumber
    m_k = \inf \Bigl\{\int_{-\infty}^{+\infty} \bigl(W(u)+(u^{(k)})^2\bigr)\, dt: u\in H^k_{\textup{loc}}(\mathbb{R}), \lim_{t\to\pm\infty} u(t)=\pm1\Bigr\}. 
\end{equation}
It is convenient to interpret this 
optimal-profile problem as an optimization on intervals and functions; to that end, we define 
\begin{eqnarray}\label{eq:def problem mk}
\widetilde m_k&:=&\inf_{T>0}\inf \Bigl\{\int_{-\infty}^{+\infty} \big(W(u)+(u^{(k)})^2\big)\, dt: u\in H^k_{\text{loc}}(\mathbb{R}),  \nonumber \\ 
    &&\hskip1.5cm u(t)=-1 \text{ for } t\leq-T, u(t)=1 \text{ for } t\geq T\Bigr\}\nonumber\\
    \nonumber&=&\inf_{T>0}m_k(T) \ =\ \lim_{T\to+\infty} m_k(T),
\end{eqnarray} 
where 
\begin{eqnarray}
   && \hskip-1cm  \nonumber m_k(T):=\inf \Bigl\{\int_{-T}^T \bigl(W(u)+(u^{(k)})^2\bigr)\, dt: u\in H^k(-T,T),   u(\pm T)=\pm 1,\nonumber \\
    && \hskip3.5cm  u^{(\ell)}(\pm T)=0 \,\,\textup{for all } \ell\in\{1,..,k-1\}\Bigr\}.\nonumber
\end{eqnarray}

As a consequence of Lemma \ref{lunghint}, we can prove that these two definitions coincide.

\begin{proposition}\label{equality} 
The equality $m_k=\widetilde m_k$ holds. 
\end{proposition}

\begin{proof}
Since all test functions for $\widetilde m_k$ are test functions for $m_k$, one inequality is trivial. 

We then consider a function $u\in H^k_{\textup{loc}}(\mathbb{R})$ such that $\lim\limits_{t\to\pm\infty} u(t)=\pm1$ and $\int_{-\infty}^{+\infty} \big(W(u)+(u^{(k)})^2\big)\, dt<+\infty$. 
With fixed $\eta$ and $N$, we can apply Lemma \ref{lunghint} twice, to the sequences $u_\e(t)= u(\frac{t-1}\e)$ and $u_\e(t)= u(\frac{t}\e)$, respectively, and with $I_{\e} = (0,1)$.  In terms of $u$, we deduce that there are $t_\e^-\in(-\frac{1}{\e}, 0)$ and $t_\e^+\in(0,\frac{1}{\e})$ such that $|u(t_\e^-)+1|<\eta$, $|u(t_\e^+)-1|<\eta$, and $|u^{(\ell)}(t_\e^\pm)|\le \frac1N$ for all $\ell\in\{1,\ldots, k-1\}$. We now choose $\eta_j\to 0$ and $N_j\to+\infty$, and consider the correspondingly defined $t_j^-$ and $t^+_j$. 

We want to extend the restriction of the test function $u$ to $(t_j^-,t_j^+)$ in order to obtain an admissible test function $u_j$ for $\widetilde m_k$ with an energy not greater than that of $u$ up to an infinitesimal term. To this end, we consider a family of auxiliary minimum problems defined as follows. 
Given $z^\ell_j$ for $\ell\in\{0,\ldots,k-1\}$, let $v_j$ be a minimizer of the problem 
\begin{eqnarray*}
\hskip.5cm \theta_j=\min\Big\{ \int_0^1(v^{(k)})^2\, dt: v^{(\ell)}(0)= z^\ell_j, v^{(\ell)}(1)=0 \hbox{ for all }\ell\in\{0,\ldots,k-1\}\Big\}, 
\end{eqnarray*}
where we used the notation $v^{(0)}=v$. 
By choosing  $z^0_j=u(t_j^+)-1$, $z^\ell_j=u^{(\ell)}(t_j^+)$ for $\ell\in\{1,\ldots,k-1\}$, we have $\theta_j\to0$, and moreover, by the uniform convergence of the corresponding minimizers $v_j$ to $0$, we also have $$\lim_{j\to+\infty}\int_0^1W(v_j(t)+1)\,dt
=0.$$ 
We then define $u_j^+$ on $(t^+_j,+\infty)$ as 
$$
u_j^+(t)=
\begin{cases} 1+v_j(t-t^+_j) &\hbox{ if $t\in(t^+_j,t^+_j+1)$},\\
1 &\hbox{ if } t\ge t^+_j+1.
\end{cases}
$$
Note that 
$$
\int_{t^+_j}^{+\infty} \big(W(u_j^+)+((u_j^+)^{(k)})^2\big)\, dt= \theta_j+\int_{0}^{1} W(v_j(t)+1)\,dt =o(1) 
$$
as $j\to+\infty$. We define symmetrically $u_j^-$ on $(-\infty, t^-_j)$, so that the functions
$$
u_j(t)=
\begin{cases} u_j^-(t) &\hbox{ if } t\le t^-_j,\\
u(t) &\hbox{ if $t^-_j<t<t^+_j$},\\
u_j^+(t) &\hbox{ if } t\ge t^+_j
\end{cases}
$$
are test functions for $\widetilde m_k$ and 
$$
\widetilde m_k\leq \int_{-\infty}^{+\infty}\big( W(u_j)+(u_j^{(k)})^2\big)\, dt\le \int_{-\infty}^{+\infty}\big( W(u)+(u^{(k)})^2\big)\, dt+ o(1)
$$
as $j\to+\infty$. 
Hence, 
$$
\widetilde m_k\leq \int_{-\infty}^{+\infty}\big( W(u)+(u^{(k)})^2\big)\, dt.
$$
By the arbitrariness of $u$ this proves the claim.
\end{proof}

\begin{proposition}\label{positivity}
We have $m_k>0$.
\end{proposition}

\begin{proof}
Since $m_k=\widetilde{m}_k$, we equivalently prove $\widetilde{m}_k>0$.

Let $\{v_j\}$ be a sequence of admissible functions for $\widetilde m_k$
such that
\begin{equation}\label{eq:vj is minimizing}\nonumber 
\lim_{j\to+\infty}\int_{-\infty}^{+\infty} \bigl(W(v_j)+(v_j^{(k)})^2\bigr)\, dt= \widetilde m_k.
\end{equation}
For each $j\geq 1$, let $T_j$ be a positive number such that $v_j(t)=-1$ for $t< -T_j$ and $v_j(t)=1$ for $t> T_j$. 

We fix $\eta<\min\{1,\sqrt{\beta_W}\})$, in view of the application of Lemma \ref{lunghint}, and consider the family of open subintervals $(a,b)$ of $(-T_j,T_j)$ such that $|v_j(a)+1|\leq \eta$, $|v_j(b)-1|\leq \eta$, and the derivatives of $v_j$ at the endpoints satisfy $|v_j^{(\ell)}|\leq 1$ for all $\ell\in\{1,\dots, k-1\}$. 
Let $(a_j,b_j)$ be a minimal such interval; that is, any subinterval of $(a_j,b_j)$ does not belong to this family. 
Hence, $(a_j,b_j)$ does not intersect the set of the points $t$ such that $||v_j(t)|-1|\leq \eta$ and $|v_j^{(\ell)}(t)|\leq 1$ for all $\ell\in\{1,\dots, k-1\}$. 

In order to apply Lemma \ref{lunghint} to prove that $b_j-a_j$ is equibounded, we use the scaling factor $\e_j=\frac{1}{2T_j}$ to define $u_j(t)= v_j(\frac{t}{\e_j}-T_j)$. 
We then have 
\[
\int_{0}^{1} \Bigl(\frac{1}{\e_j}W(u_j)+\e_j^{2k-1}(u_j^{(k)})^2\Bigr)\, dt = \int_{-T_j}^{T_j} \bigl(W(v_j)+(v_j^{(k)})^2\bigr)\, dt, 
\]
which is equibounded.  
Assume by contradiction that there exists a (not relabeled) subsequence $b_j-a_j\to+\infty$. 
Then, we can apply Lemma \ref{lunghint} to the sequence of scaled intervals $(a_j\e_j,b_j\e_j)$ 
and deduce that there exists $j_0\in\mathbb N$ such that 
$$|(a_j\e_j,b_j\e_j)\cap A_{\e_j}^\eta(1)|>0$$ 
for all $j\geq j_0$, where we recall that $A_\e^\eta(N)$ is defined as in \eqref{defaeN}. By scaling again to $(-T_j,T_j)$, this gives a contradiction.    

Up to subsequences and up to translations, we can then assume that there exists $T^\ast>0$ such that $a_j\to -T^\ast$ and  $b_j\to T^\ast$ as $j\to+\infty$. 
Up to an asymptotic negligible change of variable, it is not restrictive that $a_j=-T^\ast$ and $b_j=T^\ast$. 
By the equiboundedness of the boundary condition and by the Fundamental Theorem of Calculus we deduce that $\{v_j\}$ is equibounded in $H^k(-T^\ast,T^\ast)$.  Hence, we extract a (not relabeled) subsequence $\{v_j\}$ weakly converging to a function $v$ in $H^k(-T^\ast,T^\ast)$. 

Assume now by contradiction that $\widetilde m_k =0$. Since $\{v_j\}$ is a minimizing sequence, by the lower semicontinuity of the $L^2$-norm, we have 
\[
0 = \lim_{j\to+\infty} \int_{-T^\ast}^{T^*} \big(W(v_j)+(v_j^{(k)})^2\big)\, dt \geq \liminf_{j\to+\infty}\int_{-T^\ast}^{T^*} (v_j^{(k)})^2\, dt\geq \int_{-T^\ast}^{T^*} (v^{(k)})^2\, dt.
\]
Hence, we obtain that $v^{(k)}=0$ on $(-T^\ast,T^*)$, so that $v$ is a polynomial of degree at most $k-1$. 
Moreover, by the uniform of convergence of $v_j\to v$, we have 
\begin{equation*}
    \int_{-T^\ast}^{T^*} W(v)\, dt=\lim_{j\to+\infty}\int_{-T^\ast}^{T^\ast}W(v_j)\, dt=0, 
\end{equation*}
so that $W(v)=0$ almost everywhere, and then, since $\{W=0\}=\{-1,1\}$, we deduce that $v$ is either identically $-1$ or $1$ on $(-T^\ast,T^*)$. 
Since $|v(-T^\ast)+1|\leq \eta$ and $|v(T^\ast)-1|\leq \eta$, and $\eta<1$,  this gives a contradiction, and the proof is complete. 
\end{proof}

As it is usual in phase-transition problems, we will need some auxiliary results involving optimal profile problems approximating $m_k$ defined in  through the minimization of $\int_{I_T}(W(v)+(v^{(k)})^2)\, dt$, where $I_T$ is an interval with length $2T$ and $v$ in $H^k(I_T)$ satisfies suitable boundary conditions; in particular, we will minimize over sets  of functions with ``small'' values of the derivatives at the endpoints of the interval.   
To that end, for all $\eta\in(0,\sqrt{\beta_W})$ and for all integers $N\geq 1$, we define    
an approximating constant given by 
\begin{equation}\label{minNT}\nonumber
m_k(\eta, N):=\inf_{T>0}m_k(\eta, N,T), 
\end{equation}
where $m_k(\eta,N,T)$ is the minimum defined in \eqref{defment}.

\begin{lemma}\label{lemmaminMM} Let $\widetilde m_k$ and $m_k(\eta, N)$ be defined as above. 
Then
\begin{equation*}
     \lim\limits_{\eta\to0}\lim\limits_{N\to+\infty} m_k(\eta, N) = \widetilde m_k.
\end{equation*}
\end{lemma}

\begin{proof}

Since $m_k(\eta, N)$ is increasing as $\eta\to0$ and as $N\to+\infty$, the limits exist. 
For all $T>0$, $\eta$ and $N$, we have $m_k(T)\geq m_k(\eta, N, T)$, so that $\widetilde m_k\geq m_k(\eta, N)$ for all $\eta$ and $N$. Hence, 
$$\widetilde m_k\geq \lim_{\eta\to 0}\lim_{N\to+\infty}m_k(\eta,N).$$


We now prove the converse inequality. 
Let $\eta_j\to 0$ and $N_j\to+\infty$ be such that 
$$\lim\limits_{j\to+\infty}m_k(\eta_j,N_j)=\lim\limits_{\eta\to 0}\lim\limits_{N\to+\infty}m_k(\eta,N).$$ 
For any $j$, let $T_j>0$ and $v_j\in H^k(-T_j,T_j)$ be an admissible test function for $m_k(\eta_j,N_j,T_j)$ such that
\begin{equation*}
 \lim_{j\to+\infty} m_k(\eta_j,N_j)= \lim_{j\to+\infty}\int_{-T_j}^{T_j} \big(W(v_j)+(v_j^{(k)})^2\big)\, dt.
\end{equation*}
To define an admissible test function $u_j$ for $\widetilde m_k$, we can proceed as in the proof of Proposition \ref{equality} with $-T_j$ 
in the place of $t_j^-$, $T_j$ 
in the place of $t_j^+$, and $v_j$ in the place of the restriction of $u$ to $(t_j^-,t_j^+)$. 

In this way, the resulting $u_j$ is a test function for $\widetilde m_k$ and 
$$\lim_{j\to+\infty}m_k(\eta_j,N_j)=\lim_{j\to+\infty}\int_{-\infty}^{+\infty}\big(W(u_j)+(u_j^{(k)})^2\big)\, dt\geq \widetilde m_k,$$ 
proving the equality 
$$\lim_{\eta\to 0}\lim_{N\to+\infty}m_k(\eta,N)=\widetilde m_k,$$
and the claim. 
 \end{proof}


We introduce some quantities that will allow us to estimate oscillations of order $\eta$. 
For all $0<\eta<  \min\{1,\sqrt{\beta_W}\}$, integers $N\geq1$ and $\overline{T}>0$, we define
\begin{eqnarray}\nonumber
\hskip-0.5cm m^+_k(\eta, N, \overline T)&=&\inf_{0<T\leq \overline T}\min \Big\{\int_{0}^T \bigl(W(u)+(u^{(k)})^2\bigr)\, dt: u\in H^k(0,T), \\ 
&&\nonumber\hspace{1cm} 
|u(0)-1|\leq \eta, |u(T)-1|=2\eta,\\
&&\label{eq:def mstar +1}\hspace{1cm}\hbox{\rm and } |u^{(\ell)}(0)|\leq \frac{1}{N}\ \hbox{\rm for all } \ell\in\{1,\dots, k-1\}\Big\},
\end{eqnarray} 
and
\begin{eqnarray}\nonumber
\hskip-0.5cm m^-_k(\eta, N, \overline T)&=&\inf_{0<T\leq \overline T}\min \Big\{\int_{0}^T \bigl(W(u)+(u^{(k)})^2\bigr)\, dt: u\in H^k(0,T), \\ 
&&\nonumber\hspace{1cm} 
|u(0)+1|\leq \eta, |u(T)+1|=2\eta, \\
&&\label{eq:def mstar -1}\nonumber \hspace{1cm}\hbox{\rm and } |u^{(\ell)}(0)|\leq \frac{1}{N}\ \hbox{\rm for all } \ell\in\{1,\dots, k-1\}\Big\}.
\end{eqnarray} 

\begin{lemma} \label{lemmaminT}
    For every $\overline{T}>0$, it holds 
    \begin{equation}\label{eq:mstar +1>0} m^*_k(\eta,N,\overline{T}):=\min\{m^+_k(\eta,N,\overline{T}), m^-_k(\eta,N,\overline{T})\}>0.
    \end{equation}
\end{lemma}
\begin{proof}
    We only prove that $m^+_k(\eta,N,\overline{T})>0$, the other case being analogous. We argue by contradiction.
    
    Suppose there exists a sequence $\{T_j\}\subset(0,\overline{T}]$ and functions $v_j\in H^k(0,T_j)$ which are admissible for \eqref{eq:def mstar +1} and such that 
\begin{equation}\label{eq:minimizing mstarbis}\lim_{j\to+\infty}\int_0^{T_j}\bigl(W(v_j)+(v_j^{(k)})^2\bigr)\, dt=0.
    \end{equation}
    Without loss of generality, we may suppose that $T_j$ converges to $T\in[0,\overline{T}]$ as $j\to+\infty.$ 
    We distinguish the two cases $T=0$ and $T>0$.

\medskip
{\it Case $1$}: $(T=0)$. By \eqref{eq:minimizing mstarbis} it holds $\|v_j^{(k)}\|_{L^2(0,T_j)}\to0$ as $j\to+\infty$; hence, recalling that $v_j$ satisfies the boundary conditions appearing in \eqref{eq:def mstar +1}, we apply the Fundamental Theorem of Calculus to deduce that $\|v_j\|_{H^k(0,T_j)}\leq C$ for every $j$. We then extend each function $v_j$ to the whole interval $(0,1)$ by setting it equal to $p_j$ the polynomial of degree at most $k-1$ obtained by imposing $p_j^{(\ell)}(T_j)=v_j^{(\ell)}(T_j)$ for $\ell\in\{0,...,k-1\}$. The sequence $\{v_j\}$ is then equibounded in $H^k(0,1)$ so that there exists a subsequence, not relabeled, which converges uniformly to some $v\in H^k(0,1)$. This produces a contradiction since $T_j\to0$ and $|v_j(T_j)-v_j(0)|\geq\eta$.

\medskip

{\it Case }2: $(T>0)$. Up to an asymptotic negligible change of variable we may suppose $T_j=T$ for every $j$. Then, by \eqref{eq:minimizing mstarbis}, we have that $\{v_j\}$ is equibounded in $H^k(0,T)$ so that a certain subsequence converges uniformly on $[0,T]$ to a function $v$ in $H^k(0,T)$. Again by \eqref{eq:minimizing mstarbis} and by the lower semicontinuity of the $L^2$-norm, we have that either $v=- 1$ or $v= 1$ on $(0,T)$. But $|v_j(T_j)-1|=2\eta$ for every $j$, hence, $|v(T)-1|=2\eta$, which is a contradiction since $\eta<1$.
\end{proof}

\section{Proof of the main result}

We divide the proof of Theorem \ref{main} in proving separately compactness, lower and upper bounds.
We start by proving the equicoerciveness of the family $\{F_\e\}$.

\begin{proposition}[Compactness]\label{compactness} Let $F_\e$ be defined by \eqref{defFe} on $H^k(0,1)$, for all $\e_j\to 0$ and $\{u_j\}$ with $\sup_j F_{\e_j}(u_j)<+\infty$, there exists $u\in BV((0,1);\{-1,1\})$ such that, up to subsequences, $u_j\to u$ in measure.     
\end{proposition}

\begin{proof}
Up to subsequences, we may assume that $\sup_{\e>0} F_\e(u_\e)=S<+\infty$. We consider 
for any fixed $0<\eta<\min\{1, \beta_W\}$ and $\e>0$ the sets of indices  $\mathcal I^{1,\eta}_\e$ and 
$\mathcal I^{-1,\eta}_\e$ such that 
\begin{eqnarray*}
&&\{t\in (0,1): |u_\e(t)-1|<\eta\}=\bigcup_{j\in \mathcal I^{1,\eta}_\e}I_j,\\
&&\{t\in (0,1): |u_\e(t)+1|<\eta\}=\bigcup_{j\in \mathcal I^{-1,\eta}_\e}I_j,
\end{eqnarray*}
where each $I_j=(a_j,b_j)$ is a maximal open interval (that is, 
if $a_j\neq 0$ then $||u_\e(a_j)|-1|=\eta$, and if $b_j \neq 1 $ then $||u_\e(b_j)|-1|=\eta$, respectively) and they are pairwise disjoint. 

Let $N\in\mathbb N$ be fixed. A first remark is that for $\e$ small enough the set $(0,1)\cap A_\e^\eta(N)$ 
of the points where $||u_\e|-1|<\eta$ and $|u_\e^{(\ell)}|<\frac{1}{N\e^\ell}$ is not empty. 
This follows by the application of Lemma \ref{lunghint} to the constant sequence of intervals $I_{\e_r}=(0,1)$ and $\e_r$ any infinitesimal sequence. 
We can then define 
$$\alpha_\e=\inf\big( (0,1)\cap A_\e^\eta(N) \big)\ \ \hbox{\rm and }\ \ 
\beta_\e=\sup \big((0,1)\cap A_\e^\eta(N)\big),$$
noting that again by applying Lemma \ref{lunghint} to the sequences $\{(0,\alpha_\e)\}$ and $\{(\beta_\e,1)\}$ we have that 
\begin{equation}\label{aebe}\nonumber
|(0,1)\setminus (\alpha_\e,\beta_\e)|=o(1)_{\e\to 0}.
\end{equation} 

Let $\mathcal{I}_\e^{1,\eta}(N)$ denote the subset of indices of $\mathcal{I}_\e^{1,\eta}$ for which there exists an open interval $I\subseteq I_j=(a_j,b_j)$ such that 
\begin{equation}\label{eq:settimo}
    |u^{(\ell)}(a_j)|\leq \frac{1}{N\e^\ell} \quad \text{ and } \quad |u^{(\ell)}(b_j)|\leq \frac{1}{N\e^\ell}\quad \text{ for } \ell\in\{1,...,k-1\}.
\end{equation}
Analogously, we define $\mathcal{I}^{-1,\eta}_\e(N)$.

For each $j\in \mathcal{I}_\e^{\pm1,\eta}(N)$, we consider the maximal open interval $(a_j^{N},b_j^{N})$ for which \eqref{eq:settimo} holds true; that is, 
\begin{eqnarray*}
&\displaystyle a_j^N:=\inf\Big\{t\in[a_j, b_j) : |u^{(\ell)}(t)|\leq \frac{1}{N\e^{\ell}} \text{ for all } \ell \in\{1,...,k-1\}\Big\}, \\
&\displaystyle b_j^N:=\sup\Big\{t\in(a_j, b_j] : |u^{(\ell)}(t)|\leq \frac{1}{N\e^{\ell}} \text{ for all } \ell \in\{1,...,k-1\}\Big\}.
\end{eqnarray*}
For any $t\in (\alpha_\e,\beta_\e)\setminus \bigcup_{j\in \mathcal I_\e^{1,\eta}(N)\cup 
\mathcal I_\e^{-1,\eta}(N)}(a_j^N, b_j^N)$, we set 
\begin{eqnarray*}
&&\tau(t)=\sup\{b_j^{N}: b_j^{N}\leq t, \ j\in\mathcal I_\e^{1,\eta}(N)\cup 
\mathcal I_\e^{-1,\eta}(N)\},\\ 
&&\sigma(t)=\inf \{a_j^{N}: a_j^{N}\geq t, \ j\in\mathcal I_\e^{1,\eta}(N)\cup 
\mathcal I_\e^{-1,\eta}(N)\}. 
\end{eqnarray*}
Note that $\tau(t)$ and $\sigma(t)$ may coincide. 
For any $t$ the following properties hold: 
\begin{enumerate} 
\item[{\rm (i)}] $[\tau(t),\sigma(t)]\cap \bigcup_{j\in \mathcal I_\e^{1,\eta}(N)\cup 
 \mathcal I_\e^{-1,\eta}(N)} (a_j^N, b_j^N)=\emptyset$, so that,  recalling that $A_\e^\eta(N)$ denotes the set of points $t$ such that $||u_\e(t)|-1|< \eta$ and $|u_\e^{(\ell)}(t)|< \frac{1}{N\e^\ell}$ for every $\ell\in\{1,...,k-1\}$, it follows 
$$|(\tau(t),\sigma(t))\cap A_\e^\eta(N)|=0;$$ 
\item[\rm (ii)] if $I_j\cap [\tau(t), \sigma(t)]\neq\emptyset$ for $j\in\big(\mathcal I_\e^{1,\eta}\cup 
\mathcal I_\e^{-1,\eta}\big) \setminus \big(\mathcal I_\e^{1,\eta}(N)\cup 
\mathcal I_\e^{-1,\eta}(N)\big)$, then  
$I_j\subset [\tau(t), \sigma(t)]$; 
\item[{\rm (iii)}] if $(a_j, a_j^{N}]\cap [\tau(t), \sigma(t)]\neq\emptyset$, then  
$(a_j, a_j^{N}]\subset[\tau(t), \sigma(t)]$; the same property holds for $[b_j^{N}, b_j)$; 
\item[{\rm (iv)}] by the continuity of $u_\e$ and of its derivatives,  
$$||u_\e(\tau(t))|-1|\leq \eta \ \ \hbox{\rm and } \ \ |u_\e^{(\ell)}(\tau(t))|\leq \frac{1}{N\e^\ell}  \ \ \hbox{\rm for all }\ \ell\in\{1,\dots, k-1\},$$ 
and the same properties hold in $\sigma(t)$.  
\end{enumerate}  
Moreover, for different $t$ the intervals $(\tau(t), \sigma(t))$ either coincide or are disjoint; indeed, 
if $s\in (\tau(t), \sigma(t))$ for some $t$, then $(\tau(t), \sigma(t))=(\tau(s), \sigma(s))$.

Then, we can write 
the set $(a_\e,b_\e)\setminus   \bigcup_{j\in \mathcal I_\e^{1,\eta}(N)\cup 
\mathcal I_\e^{-1,\eta}(N)}(a_j^N, b_j^N)$ as a disjoint union of intervals of the form $[\tau(t), \sigma(t)]$; that is, there exists a set of indices $\Lambda_\e(\eta,N)$ such that
$$(a_\e,b_\e)\setminus   \bigcup_{j\in \mathcal I_\e^{1,\eta}(N)\cup 
\mathcal I_\e^{-1,\eta}(N)}(a_j^N, b_j^N)=\bigcup_{\lambda\in \Lambda_\e(\eta, N)}[\tau_\lambda, \sigma_\lambda],$$ 
where the interior of intervals are pairwise disjoint and the endpoints $\tau_\lambda$ and $\sigma_\lambda$
belong to the sets $\{\tau(t): t\in (a_\e,b_\e)\setminus   \bigcup_{j\in \mathcal I_\e^{1,\eta}(N)\cup 
\mathcal I_\e^{-1,\eta}(N)}(a_j^N,b_j^N)\}$ and $\{\sigma(t): t\in (a_\e,b_\e)\setminus   \bigcup_{j\in \mathcal I_\e^{1,\eta}(N)\cup 
\mathcal I_\e^{-1,\eta}(N)}(a_j^N,b_j^N)\}$, respectively. 

Since the measure of the set of points $t$ such that $||u_\e(t)|-1|\geq\eta$ is bounded by $\e \frac{S}{\eta^2}$, we deduce that there exists a countable subfamily $\Lambda^\ast_\e(\eta, N)\subset \Lambda_\e(\eta, N)$ such that 
\begin{equation}\label{comp}\nonumber
(a_\e,b_\e)\setminus   \bigcup_{j\in \mathcal I_\e^{1,\eta}(N)\cup 
\mathcal I_\e^{-1,\eta}(N)}(a_j^N, b_j^N)=J_\e\cup \bigcup_{\lambda\in \Lambda^\ast_\e(\eta, N)}(\tau_\lambda, \sigma_\lambda),
\end{equation} 
where $|J_\e|\to 0$ as $\e\to 0$.  

We now prove that there exists a constant $T_k(\eta, N, S)>0$ 
such that 
\begin{equation}\label{stimatlsl}
|(\tau_\lambda, \sigma_\lambda)|\leq T_k(\eta, N, S)\e
\end{equation}  
for any $\e>0$ and $\lambda\in\Lambda^\ast_\e(\eta, N)$. 
  
We reason by contradiction. If the claim does not hold, then for any $r\in\mathbb N$
there exist $\e_r>0$ and $\lambda_r\in \Lambda^\ast_{\e_r}(\eta, N)$ such that 
$$|(\tau_{\lambda_r}, \sigma_{\lambda_r})|> r\e_r.$$ 
Since $(\tau_{\lambda_r}, \sigma_{\lambda_r})\subset (0,1)$, this implies that $\e_r\to 0$ as $r\to+\infty$. 
We can then apply Lemma \ref{lunghint} 
to the family of intervals $I_{\e_r}=(\tau_{\lambda_r}, \sigma_{\lambda_r}),$ deducing that 
there exists $r_0\in\mathbb N$ such that 
$$|A_{\e_r}^\eta(N)\cap (\tau_{\lambda_r}, \sigma_{\lambda_r})|>0$$
for any $r\geq r_0$. This gives a contradiction, since property (i) holds for any $(\tau_{\lambda_r}, \sigma_{\lambda_r})$.

We say that an interval $(\tau_\lambda, \sigma_\lambda)$ with $\lambda\in \Lambda^\ast_{\e}(\eta, N)$ is an {\em effective transition interval} if 
either $|u_\e(\tau_\lambda)-1|\leq \eta$ and $|u_\e(\sigma_\lambda)+1|\leq \eta$, or 
$|u_\e(\tau_\lambda)+1|\leq \eta$ and $|u_\e(\sigma_\lambda)-1|\leq \eta$. We now prove that the number of such intervals is equibounded.  Indeed, setting $v_\e(s)=u_\e(\e s+\frac{\sigma_\lambda+\tau_\lambda}{2})$, by a change of variable we obtain
\begin{equation}\label{boundtransition}
F_\e(u_\e;(\tau_\lambda, \sigma_\lambda))=\int_{-\frac{\sigma_\lambda-\tau_\lambda}{2\e}}^{\frac{\sigma_\lambda-\tau_\lambda}{2\e}}
\Big(W(v_\e) +(v_\e^{(k)})^2\Big)\, ds
\geq m_k(\eta, N)\ \geq \ \frac{m_k}{2},
\end{equation} 
where the last inequality follows by Lemma \ref{lemmaminMM}. Since $m_k>0$, the number of intervals of effective transitions is equibounded by a constant which is independent on $\e, \eta$ and $N$; more precisely, by $\frac{2S}{m_k}$.

We now consider intervals $(\tau_\lambda, \sigma_\lambda)$ which are not intervals of effective transitions; that is, either 
$|u_\e(\tau_\lambda)-1|\leq \eta$ and $|u_\e(\sigma_\lambda)-1|\leq \eta$ or $|u_\e(\tau_\lambda)+1|\leq \eta$ and $|u_\e(\sigma_\lambda)+1|\leq \eta$. 
Let $(\tau_\lambda,\sigma_\lambda)$ be such that 
$|u_\e(\tau_\lambda)-1|\leq \eta$ and $|u_\e(\sigma_\lambda)-1|\leq \eta$, and suppose that there exists a point $t$ in $(\tau_\lambda,\sigma_\lambda)$ such that $|u_\e(t)-1|\geq 2\eta$. 
Then, setting $t^\prime=\inf \{t:|u_\e(t)-1|\geq 2\eta\}$, by a change of variable we have that 
\begin{eqnarray*}
F_\e(u_\e;(\tau_\lambda, \sigma_\lambda))&\geq&F_\e(u_\e;(\tau_\lambda, t^\prime))
=\int_{\tau_\lambda}^{t^\prime}
\Big(\frac{1}{\e}W(u_\e) +\e^{2k-1}(u_\e^{(k)})^2\Big)\, dt\\
&=&\int_{0}^{\frac{t^\prime-\tau_\lambda}{\e}}
\Big(W(v_\e) +(v_\e^{(k)})^2\Big)\, ds\\
&\geq&\inf_{T\leq T_k(\eta, N, S)} \min_v\Big\{\int_{0}^{T}
\Big(W(v_\e) +(v_\e^{(k)})^2\Big)\, ds\Big\},
\end{eqnarray*} 
where the minimum is taken over $v\in H^{k}(0,T)$ such that $|v(T)-1|=2\eta$, $|v(0)-1|\geq \eta$, $v^{(\ell)}(0)\leq \frac{1}{N}$ for all $\ell\in\{1,\dots, k-1\}$, and we used the fact that $|(\tau_\lambda,\sigma_\lambda)|\leq T_k(\eta, N, S)\e$ by \eqref{stimatlsl}. 
Then, by applying Lemma \ref{lemmaminT} with $\eta>0$ and $\overline{T}=T_k(\eta, N,S)$, we obtain 
\begin{equation*}
F_\e(u_\e;(\tau_\lambda, \sigma_\lambda))\geq m^\ast_k(\eta,N,T_k(\eta, N,S))>0. 
\end{equation*} 
Since a corresponding argument holds if $|u_\e(\tau_\lambda)+1|\leq \eta$ and $|u_\e(\sigma_\lambda)+1|\leq \eta$, we can conclude that 
\begin{equation}\label{prox}\nonumber
||u_\e(t)|-1|\leq 2\eta \ \ \hbox{\rm in } \ (\tau_\lambda,\sigma_\lambda) 
\end{equation}
except for a finite number of indices $\lambda$ (which of course depends on $\eta$ and $N$)  We let $c_k(\eta, N)$ denote the number of such $\lambda$.  

Let now $\{(\tau_i,\sigma_i)\}_{i=1}^{M_\e}$ with $\tau_i<\tau_{i+1}$ denote the set of the effective transition intervals. Since $M_\e\leq \frac{2S}{m_k}$, we may suppose $M_\e$ is independent of $\e$, up to passing to a subsequence, and we set $M_\e=M$. 
We construct a sequence of functions $u_\e^{\eta, N}$ in $(0,1)$ as follows.  
To give a definition in the intervals $(\sigma_i,\tau_{i+1})$, we note that if there exists 
an index $j\in \mathcal I_\e^{-1,\eta}(N)$ such that 
$(a_j^N, b_j^N)\subset (\sigma_i,\tau_{i+1})$, then $(\sigma_i,\tau_{i+1})\cap\bigcup_{\mathcal I_\e^{1,\eta}(N)}(a_j^N, b_j^N)=\emptyset$. Indeed, otherwise $(\sigma_i,\tau_{i+1})$ would contain an effective transition interval. 
We then set 
$$u_\e^{\eta, N}=-1 \ \hbox{\rm in } (\sigma_i,\tau_{i+1})$$ 
if there exists 
$j\in \mathcal I_\e^{-1,\eta}(N)$ such that 
$(a_j^N, b_j^N)\subset (\sigma_i,\tau_{i+1})$, and 
$$u_\e^{\eta, N}=1$$ 
otherwise in $(a_\e,b_\e)$. Finally, we just extend by continuity on the interval $(0,1)$. Note that, with this construction, it holds
\begin{equation}\label{boundjump}
\#S(u_\e^{\eta,N}) = M, \quad \text{ for every } \e, \eta, N.
\end{equation} 
Moreover, by \eqref{boundtransition} we have that 
\begin{equation}\label{liminfueta}
\liminf_{\e\to 0}F_\e(u_\e)\geq M m_k(\eta, N). 
\end{equation}
Up to subsequences, by \eqref{boundjump} we can suppose that $u_\e^{\eta_\e, N_\e}\to u$ in $L^1$,  the limit $u$ takes values in $\{-1,1\}$, and its jump set $S(u)$ is composed of not more than $M$ points. 

By construction, we also have that 
\begin{equation} \label{finalcompactness1}
|\{t\in(0,1): |u_\e(t)-u_\e^{\eta, N}(t)|>2\eta\}| \leq  C_k(\eta,N)\,\e,
\end{equation}
where $C_k(\eta,N):=\frac{2S}{m_k}+ c_k(\eta,N)$. 
By this estimate, we show that, up to subsequences, $\{u_\e\}$ converges in measure to the same limit $u$. 
Indeed, with fixed $\delta>0$, by the convergence of $\{u_\e^{\eta_\e,N_\e}\}$ we have that
\begin{equation}\label{finalcompactness3}\nonumber
|\{|u_\e^{\eta_\e,N_\e}-u|>\delta\}|=o(1) 
\end{equation}
as $\e\to0$.
By a triangular inequality and \eqref{finalcompactness1}, for $\e$ small enough  
we have
\begin{eqnarray*}
    |\{|u_\e-u|>2\delta\}| &\leq&  |\{|u_\e-u_\e^{\eta_\e,N_\e}|>2\eta_\e\}| + |\{|u_\e^{\eta_\e,N_\e}-u|>\delta\}| \\
    &\le& C_k(\eta_\e, N_\e) \e+o(1)=o(1)
\end{eqnarray*}
as $\e\to 0$, concluding the proof. 
\end{proof}

Finally, we prove the $\Gamma$-convergence result.

\begin{proposition}[Lower bound] \label{prop:gammaliminf1d} Let $u\in BV((0,1);\{-1,1\})$ and $u_\e\to u$ in measure. Then
    \begin{equation*}
        \liminf_{\e\to 0}F_\e(u_\e) \geq m_k \# S(u).
    \end{equation*}
    
\end{proposition}
\begin{proof}
    From the last part of the proof and by \eqref{liminfueta} we deduce that, up to the extraction of subsequences,
\begin{equation}\label{limN}\nonumber
\liminf_{\e\to 0} F_\e(u_\e)\ge m_k(\eta,N) \#S(u).
\end{equation}
We can then let $N\to+\infty$ and $\eta\to0$ and apply Lemma \ref{lemmaminMM} and Proposition~\ref{equality}. 
\end{proof}

\begin{proposition}[Upper bound]
Let $u\in BV((0,1);\{-1,1\})$. Then there exists a recovery sequence $u_\e\to u$ in $L^1(0,1)$ such that
    \begin{equation*}
        \limsup_{\e\to 0^+}F_\e(u_\e)=m_k \# S(u).
    \end{equation*}
\end{proposition}
\begin{proof}
    Let $u\in BV((0,1);\{-1,1\})$ and denote by $x_1,...,x_r$ its jump points. We set $J_1:=\{n\in\{1,...,r\} : u(x_n^{+})-u(x_n^{-})=2 \}$, $J_2:=\{n\in\{1,...,r\} : u(x_n^{+})-u(x_n^{-})=-2 \}$, and for convenience, we also set $x_0=0$ and $x_{r+1}=1$. For a fixed real positive $\delta<\min\{x_{n}-x_{n-1}:\,\, n=1,...,r+1\}$, consider a function $v$ admissible for $\widetilde m_k$ such that
    \begin{equation}\label{eq:minimality in limsup}
        \int_{-\infty}^{+\infty} \big(W(v)+(v^{(k)})^2\big) dt\leq m_k+\frac{\delta}{r}.
    \end{equation}
    Recall that there exists  $T>0$ such that $v=-1$ in $(-\infty,-T)$ and $v=1$ in $(T,+\infty)$. 
    
    Let $\e_j \to0^+$ and suppose that $\frac{\delta}{2\e_j}>T$ for every $j$. For $n\in\{1,..,r\}$, define intervals $I_n:=[(x_{n-1}+x_n)/2,(x_n+x_{n+1})/2]$ and consider the sequence defined as 
    \begin{equation*}
       u_j(t):= \begin{cases}
        \displaystyle v\Bigl(\frac{t-x_n}{\e_j}\Bigr)\quad &\text{if } t\in I_n\text{ and } n\in J_1,\\[7pt]
        \displaystyle v\Bigl(-\frac{t-x_n}{\e_j}\Bigr)\quad & \text{if } t\in I_n \text{ and }n\in J_2,\\[7pt]
        u(t)\quad& \text{otherwise}.
        \end{cases}
    \end{equation*}
\noindent We note that the condition $\frac{\delta}{2\e_j}>T$ implies that $\{u_j\}\subset H^k(0,1)$; and moreover, $u_j\to u$ in $L^1(0,1)$ as $j\to+\infty$. To conclude, we estimate
    \begin{eqnarray*}
        \lim_{j\to+\infty}F_{\e_j}(u_j)&=&\lim_{j\to+\infty}\sum_{n=1}^r\int_{I_n}
        \Big(\frac{1}{\e_j}W(u_j(t))+\e_j^{2k-1}(u_j^{(k)}(t))^2\Big) dt\\
        &=&\lim_{j\to+\infty}\sum_{n\in J_1}\!\int_{\frac{x_{n-1}-x_n}{2\e_j}}^{\frac{x_{n+1}-x_n}{2\e_j}} \big(W(v(s))+(v^{(k)}(s))^2\big) ds \\
        && + \lim_{j\to+\infty}\sum_{n\in J_2}\int_{\frac{x_{n-1}-x_{n}}{2\e_j}}^{\frac{x_{n+1}-x_{n}}{2\e_j}}\big( W(v(-s))+(v^{(k)}(-s))^2\big) ds\\
        &=&\sum_{n\in J_1}\int_{-\infty}^{+\infty} \big(W(v(s))+(v^{(k)}(s))^2\big) ds\\
        &&+\sum_{n\in J_2}\int_{-\infty}^{+\infty} \big(W(v(-s))+(v^{(k)}(-s))^2\big) ds\\
        &\leq& m_k\# S(u)+\delta,
    \end{eqnarray*}
   where in the last inequality we have used \eqref{eq:minimality in limsup}. The arbitrariness of $\delta$ and a diagonal argument allow us to conclude.
\end{proof}

\section{$\Gamma$-limit in dimension $d>1$}
%

{
In this section we prove Theorem \ref{main-d} in the general $d$-dimensional case. 
To that end, we use some standard notions in the theory of sets of finite perimeter, for which we refer to \cite{Maggi} (see also \cite{BLN98}).
Moreover, we use some notation for $d$-dimensional sets.
We let $Q_1=(-\frac12,\frac12)^{d }$, and more in general, we let $Q_\rho^\nu(x)$ denote a cube with centre $x$, side length $\rho$ and one face orthogonal to $\nu$. In particular we may take $Q^{e_j}_\rho(0)=Q_\rho$ if  $x=0$ and $\nu =e_j$, the $j$-th element of the canonical basis of $\R^d$.} 

%

\medskip
{
We first prove the lower bound using the blow-up technique by Fonseca and M\"uller \cite{FM} (see also \cite{BMS}) combined with a slicing argument that allows us to reduce to the one-dimensional case in Theorem \ref{main}.}

\begin{proposition} [Lower bound]
    Let $u\in BV(\Omega;\{-1,1\})$ and let $\{u_\e\}$  be a sequence in $H^k(\Omega)$ such that $u_\e\to u$ in $L^1(\Omega)$. Then
\begin{equation*}
    \liminf_{\e\to 0}F_\e(u_\e)\geq m_k \mathcal{H}^{d-1}(S(u)).
\end{equation*}
\end{proposition}

\begin{proof}
{
    Without loss of generality, we assume $\sup_\e F_\e(u_\e)<+\infty$. Let $\e_j\to0^+$ be such that $$\lim_{j\to+\infty}F_{\e_j}(u_{\e_j})=\liminf_{\e\to 0}F_\e(u_\e)$$ and let $\{u_{\e_j}\}$ be a sequence converging to $u$ in $L^1(\Omega)$. Since $u\in BV(\Omega;\{-1,1\})$, the set $E=\{x\in\Omega: u(x)=1\}$ is a set of finite perimeter in $\Omega$ such that $u=2\chi_E-1$ and $S(u)=\Omega\cap\partial^*E$, where $\partial^*E$ is the reduced boundary of $E$, up to an $\mathcal{H}^{d-1}$-negligible set. This is a $(d-1)$-rectifiable set, and for every $x\in\partial^*E$ there exists $\nu_E(x)$, the inner unit normal in the measure-theoretical sense.

    For every positive integer $j$, we introduce the measure $\mu_{\e_j}$ defined for $B\subset \Omega$ Borel as
    \begin{equation*}\label{eq:Definition of measures muepsilon}
        \mu_{\e_j}(B):= \int_B\Big(\frac{1}{\e_j}W(u_{\e_j})+\e_j^{2k-1}\|D^ku_{\e_j}\|^2\Big) dx.
    \end{equation*}
    Since $\sup_j \mu_{\e_j}(\Omega)<+\infty$, up to a not relabeled subsequence, we may assume that $\{\mu_{\e_j}\}$ weakly* converges to a finite Borel measure $\mu$. 
By the Besicovitch Derivation Theorem, for $\mathcal{H}^{d-1}$-almost every $x$ in $\partial^*E$, it holds
    \begin{equation}\label{eq:Besicovitch mu}
        \frac{d\mu}{d(\mathcal{H}^{d-1}\mres \partial^*E)}(x)=\lim_{\rho\to 0}\frac{\mu(Q_\rho^{\nu}(x))}{\mathcal{H}^{d-1}(Q_\rho^\nu(x))\cap \partial^*E)}=\lim_{\rho\to 0}\frac{\mu(Q_\rho^\nu(x))}{\rho^{d-1}},
    \end{equation}
where $\nu=\nu_E(x)$ is the inner unit normal to $E$ at $x$ and in the last equality we have used that the $(d-1)$-dimensional density of $\partial^*E$ is equal to $1$ for $\mathcal{H}^{d-1}$-almost every point in the reduced boundary of $E$.
Since $\mu$ is a finite measure, we also have
\begin{equation}\label{eq: is at most countable}
    \lim_{j\to +\infty}\mu_{\e_j}(Q_\rho^\nu(x)) = \mu(Q_\rho^\nu(x))
    \end{equation}
for all $\rho>0$ but a countable set. Hence, if we fix a point $x\in \partial^*E$ such that \eqref{eq:Besicovitch mu} holds, a diagonal argument and \eqref{eq: is at most countable} imply that there exists $\rho_j=\rho(\e_j)\to0^+$ such that
\begin{equation}\label{eq:cosa da calcolare}
        \frac{d\mu}{d(\mathcal{H}^{d-1}\mres \partial^*E)}(x)=\lim_{j\to+\infty}\frac{\mu_{\e_j}(Q_{\rho_j}^\nu(x))}{\rho_j^{d-1}}
    \end{equation}
and $\frac{\e_j}{\rho_j}\to0$.

Up to a translation and rotation argument, we may suppose that $x=0\in \partial^*E$ is such that \eqref{eq:cosa da calcolare} holds and that $\nu_E(0)=e_d$. In this case, the following blow-up condition holds:
     \begin{equation}\label{eq:blow-up}
         \frac{1}{\rho_j}\chi_E\xrightarrow{L^1_{\text{loc}}(\R^d)} \chi_{\{x_d\geq0\}}=:H \,\,\,\text{as } j\to+\infty,
     \end{equation}
     where $H$ is the half-space given by $x_d\geq 0$.
  
Applying the change of variable $y:=x/\rho_j$ and setting $v_{\e_j}(y):=u_{\e_j}(x)$, we get
    \begin{eqnarray*}
        \frac{\mu_{\e_j}(Q_{\rho_j})}{\rho_j^{d-1}}&=&\frac{1}{\rho_j^{d-1}}\int_{Q_{\rho_j}}\Big(\frac{1}{\e_j}W(u_{\e_j}(x))+\e_j^{2k-1}\|D^k u_{\e_j}(x)\|^2\Big)dx\\   &=&\int_{Q_1}\Big(\frac{1}{\delta_j}W(v_{\e_j}(y))+\delta_j^{2k-1}\|D^k v_{\e_j}(y)\|^2\Big)dy,
    \end{eqnarray*}
    where $\delta_{j}:=\frac{\e_j}{\rho_j}$. We set $y=(\xi,t)$, with $\xi\in [-1/2,1/2]^{d-1},\, t\in[-\frac{1}{2},\frac{1}{2}]$ and $v_{\e_j}^{\xi}(t):=v_{\e_j}(\xi,t)$, and we note that
    \[
    \|D^kv_{\e_j}\|\geq \Bigl|\frac{d^k}{dt^k}v_{\e_j}^\xi(t)\Bigr|.
    \]
    Therefore, we have
    \begin{equation*}
        \frac{\mu_{\e_j}(Q_{\rho_j})}{\rho_j^{d-1}} \geq \int_{[-\frac{1}{2},\frac{1}{2}]^{d-1}}\int_{-\frac{1}{2}}^{\frac{1}{2}}\Big(\frac{1}{\delta_j} W(v^\xi_{\e_j}(t)) +\delta_j^{2k-1} \Bigl(\frac{d^k}{dt^k}v_{\e_j}^\xi(t)\Bigr)^2\Big) dt\,d\xi
    \end{equation*}
and by Fatou's Lemma 
\begin{eqnarray*} \nonumber
 &&\liminf_{j\to+\infty} \frac{\mu_{\e_j}(Q_{\rho_j})}{\rho_j^{d-1}} \\  
&\geq& \int_{[-\frac{1}{2},\frac{1}{2}]^{d-1}}\liminf_{j\to+\infty}\int_{-\frac{1}{2}}^{\frac{1}{2}}\Big(\frac{1}{\delta_j} W(v^\xi_{\e_j}(t)) +\delta_j^{2k-1} \Bigl(\frac{d^k}{dt^k}v_{\e_j}^\xi(t)\Bigr)^2\Big) dt\,d\xi. 
\end{eqnarray*}
Since $\{u_{\e_j}\}$ converges to the function $2\chi_E-1$ in $L^1(\Omega)$, by \eqref{eq:blow-up} we have that $v_{\e_j}\to 2H-1$ in $L^1(Q_1)$, and in turn, $v^\xi_{\e_j} \to 2\chi_{[0,\frac{1}{2})}-1=:v^\xi$ for Lebesgue-almost every $\xi\in [-\frac{1}{2},\frac{1}{2}]^{d-1}$. We apply the lower estimate in Theorem \ref{main} to infer that 
\[
\liminf_{j\to+\infty} \frac{\mu_{\e_j}(Q_{\rho_j})}{\rho_j^{d-1}} \geq \int_{[-\frac{1}{2},\frac{1}{2}]^{d-1}} m_k \# S(v^\xi) \,d\xi=m_k.
\]

Since this argument is valid for all $x\in \partial^*E$, recalling \eqref{eq:cosa da calcolare}, we obtain
\begin{eqnarray*}
        \mu(\Omega) &\geq& \int_{\partial^*E} \frac{d\mu}{d(\mathcal{H}^{d-1}\mres \partial^*E)}(x) d\mathcal{H}^{d-1}(x) \\
        &\geq& m_k \mathcal{H}^{d-1}(\partial^*E) = m_k \mathcal{H}^{d-1}(S(u)),
\end{eqnarray*}
and since the weak$^*$ convergence of $\mu_{\e_j}$ to $\mu$ implies that $\liminf\limits_{j\to+\infty}\mu_{\e_j}(\Omega) \geq \mu(\Omega)$, finally we get 
\begin{equation*}
    \liminf_{\e\to 0} \mu_\e(\Omega) =\lim_{j\to+\infty} \mu_{\e_j}(\Omega) \geq \mu(\Omega) \\
    \geq m_k\mathcal{H}^{d-1}(S(u)),
\end{equation*}
and the proof is complete.}\end{proof}

We conclude the proof of the $\Gamma$-convergence by showing the existence of a recovery sequence. In what follows $C$ denotes a generic positive constant which may vary from line to line.

\begin{proposition}[Upper bound]
    Let $u\in BV(\Omega;\{-1,1\})$. Then there exists a recovery sequence $u_\e\in H^k(\Omega)$ converging to $u$ in $L^1(\Omega)$ such that
    \begin{equation}\label{eq:}
        \limsup_{\e\to0}F_\e(u_\e)\leq m_k \mathcal{H}^{d-1}(S(u)).
    \end{equation}
\end{proposition}
\begin{proof}
Let $E$ be the set of finite perimeter in $\Omega$ such that $u=2\chi_E-1$. Then, there exists a sequence $\{E_j\}$ of  bounded open sets in $\mathbb{R}^d$ with $C^\infty$ boundary such that $\chi_{E_j\cap\Omega}\to\chi_{E}$ in $L^1(\Omega)$ and $\textup{Per}(E_j;\Omega)\to \textup{Per}(E;\Omega)$. Therefore, since $\mathcal{H}^{d-1}(S(u))=\textup{Per}(E;\Omega)$, it is enough to prove the $\Gamma$-limsup inequality \eqref{eq:} only for functions $u=2\chi_E-1$ when $E=\widetilde E \cap \Omega$ for some smooth bounded set  $\widetilde E$ in $\mathbb{R}^d$. We set $M :=\partial \widetilde{E}$.

Let $\e_j\to0^+$, and consider $\{v_j\}$ a sequence of functions in $H^k_{\rm loc}(\R)$ such that 
\begin{equation}\label{eq:vj is an optimal profile limd}
    \lim_{j\to+\infty}\int_{-\infty}^{+\infty}\Big( W(v_j)+(v_j^{(k)})^2\Big)dt=m_k
\end{equation}
and
\begin{equation}\label{eq:def vj limd}\nonumber
     v_j(t)= \begin{cases}
      -1 & \text{if }t\leq - \frac{1}{\sqrt{\e_j}}, \\[5pt] 
       1  & \text{if }t\geq \frac{1}{\sqrt{\e_j}}.
    \end{cases}
\end{equation}
Since $M$ is a smooth manifold in $\mathbb{R}^d$, there exists $\delta_0>0$ such that for every $\delta<\delta_0$ there exists a unique smooth projection from $U_\delta$ on $M$, where $U_\delta:=\{x\in\R^d:\, d(x,M)< \delta \}$ is a tubular neighborhood of $M$. Therefore, denoting by $\widetilde{d}_{M}$ the signed distance function from $M$ which is negative inside of $E$, we have that $\widetilde d_M\in C^k(U_\delta)$ for every $\delta<\delta_0$.

We suppose that $\sqrt{\e_j}<\delta_0/2$ for every $j$ and define  
\begin{equation}\label{vlimsup}
   u_j(x)= \begin{cases}
   -1  & \text{if } x\in\Omega \setminus(E\cup U_j),\\[5pt]
    v_j\bigl(\frac{\widetilde{d}_M(x)}{\e_j}\bigr) & \text{if }x\in \Omega\cap U_j,\\[5pt]
    1  & \text{if } x\in E\setminus U_j,\\        
    \end{cases}
\end{equation}
where $U_j:=U_{\sqrt{\e_j}}$.

In order to estimate $F_{\e_j}(u_{\e_j})$, we compute the derivatives of order $k$ of the function $u_j$. Let $\alpha$ be a multi-index of order $k$. We observe that
\begin{eqnarray*} \nonumber
 &&\hskip-.5cm \partial^\alpha u_j =  \frac{1}{\e^k}v_j^{(k)}\Bigl(\frac{\widetilde d_M}{\e}\Bigr) (\partial_{\alpha_1} \widetilde d_M )\cdots (\partial_{\alpha_k} \widetilde d_M )  + \sum_{\ell=1}^{k-1} \frac{1}{\e^\ell} v_j^{(\ell)}\Bigl(\frac{\widetilde d_M}{\e}\Bigr) P^\ell_\alpha(A(\partial^{\alpha} \widetilde d_M)),
\end{eqnarray*}
where $A(\partial^\alpha \widetilde d_M)$ represents the ordered vector of all the derivatives $\partial^\beta \widetilde d_M$ with $\beta \subset\alpha, |\beta|\geq1$, and $P^\ell_\alpha : \R^{2^{|\alpha|}-1} \to \R$ are polynomials of degree at most $k$. Therefore, the $k$-tensor $D^ku_j$ can be written in the form
\begin{equation}\label{eq:tensore dk}
    D^ku_j =  \frac{1}{\e^k}v_j^{(k)}\Bigl(\frac{\widetilde d_M}{\e}\Bigr)(\nabla \widetilde d_M \otimes ... \otimes\nabla \widetilde d_M) + \sum_{\ell=1}^{k-1} \frac{1}{\e^\ell} v_j^{(\ell)}\Bigl(\frac{\widetilde d_M}{\e}\Bigr) T^\ell,
\end{equation}
where $T^\ell,\, \ell\in\{1,...,k-1\}$ are $k$-tensors having equibounded operatorial norm. Indeed, the coefficients of $T^\ell$ only depend on the derivatives of $\widetilde d_M$ up to order $k$, and these are uniformly bounded on $U_j$ since we assumed $\sqrt{\e_j}< \delta_0/2$, which implies the existence of a compact set $K$ such that $U_j \subset \subset K \subset U_{\delta_0}$.

By \eqref{vlimsup} and \eqref{eq:tensore dk} we obtain 
\begin{eqnarray*}
   \nonumber && \hskip-0.8cm \limsup_{j\to+\infty}F_{\e_j}(u_j)
   =\limsup_{j\to+\infty} \bigg(\int_{U_j}\frac{1}{\e_j}W\Bigl(v_j\Bigl(\frac{\widetilde d_M}{\e_j}\Bigr)\Bigr)\,dx\\    &&\hskip-.4cm +\e_j^{2k-1}\int_{U_j}\Bigl\|\frac{1}{\e_j^k}v_j^{(k)}\Bigl(\frac{\widetilde d_M}{\e_j}\Bigr)(\nabla \widetilde d_M \otimes ... \otimes\nabla \widetilde d_M)  + \sum_{\ell=1}^{k-1} \frac{1}{\e_j^\ell} v_j^{(\ell)}\Bigl(\frac{\widetilde d_M}{\e_j}\Bigr) T^\ell\Bigr\|^2 dx\bigg).
\end{eqnarray*}
Let $\Phi :M\times(-\delta_0/2,\delta_0/2)\to U_{\delta_0/2}$ be the diffeomorphism defined by $\Phi(y,t):=y+t\nu(y)$, where $\nu(y)$ is the outer unit normal to $M$ at $y$, and denote by $J(y,t)$ its Jacobian. Applying this change of variable, we get
\begin{eqnarray*}
    \nonumber &&\hskip-1cm\limsup_{j\to+\infty}F_{\e_j}(u_j) = \limsup_{j\to+\infty}\int_{M\cap \Omega}\int_{-\sqrt{\e_j}}^{\sqrt{\e_j}}\bigg(\frac{1}{\e_j}W\Bigl(v_j\Bigl(\frac{t}{\e_j}\Bigr)\Bigr) \\
    && \hskip2cm +\e_j^{2k-1}\Bigl\|\frac{1}{\e_j^k}v_j^{(k)}\Bigl(\frac{t}{\e_j}\Bigr)(\nabla \widetilde d_M \otimes \cdots \otimes\nabla \widetilde d_M)(\Phi(y,t))  \\
    &&\hskip2cm+ \sum_{\ell=1}^{k-1} \frac{1}{\e_j^\ell} v_j^{(\ell)}\Bigl(\frac{t}{\e_j}\Bigr) T^\ell(\Phi(y,t))\Big\|^2\bigg) J(y,t)\, dtd\mathcal{H}^{d-1}(y),
\end{eqnarray*}
and taking into account that the Jacobian and the operatorial norm of all $T^\ell$ are bounded on $M\times (-\delta_0/2,\delta_0/2)$, we have
$$
\limsup_{j\to+\infty}F_{\e_j}(u_j)\le\limsup_{j\to+\infty} \bigg(I^{k}_j+C\sum_{\ell=1}^{k-1}I_j^{\ell}+C\sum_{\substack{{r,\ell=1}\\r< \ell}}^kI^{{r,\ell}}_j\bigg), 
$$
where
\begin{eqnarray*}
 &&   I^{k}_j:=\int_{M\cap\Omega}\int_{-\sqrt{\e_j}}^{\sqrt{\e_j}}
    \Big(\frac{1}{\e_j}W\Bigl(v_j\Bigl(\frac{t}{\e_j}\Bigr)\Bigr)+\frac{1}{\e_j}\Bigl(v_j^{(k)}\Bigl(\frac{t}{\e_j}\Bigr)\Bigr)^2\Big)J(y,t)\,dt\,d\mathcal{H}^{d-1}(y)
    \\&& I_j^{\ell}:=\int_{-\sqrt{\e_j}}^{\sqrt{\e_j}}\e_j^{2k-1-2\ell}\Bigl(v_j^{(\ell)}\Bigl(\frac{t}{\e_j}\Bigr)\Bigr)^2\,dt
     \\
  &&  I^{{r,\ell}}_j:=\int_{-\sqrt{\e_j}}^{\sqrt{\e_j}}\e_j^{2k-1-r-\ell}\Bigl|v_j^{(r)}\Bigl(\frac{t}{\e_j}\Bigr)v_j^{(\ell)}\Bigl(\frac{t}{\e_j}\Bigr)\Bigr|dt,
 \end{eqnarray*}
and we have also used that $\|\nabla \widetilde d_M \otimes\cdots\otimes\nabla \widetilde d_M \|\leq |\nabla \widetilde d_M |^k=1$ since $\widetilde{d}_M$ is a $1$-Lipschitz function. 

We first prove that
\begin{equation*}
    \limsup_{j\to+\infty} I_j^k \leq m_k \mathcal{H}^{d-1}(S(u)).
\end{equation*}
By a change of variable, we have
\begin{eqnarray*}
    I_j^k & =&  \int_{M\cap\Omega}\int_{-\sqrt{\e_j}}^{\sqrt{\e_j}}
    \Big(\frac{1}{\e_j}W\Bigl(v_j\Bigl(\frac{t}{\e_j}\Bigr)\Bigr)+\frac{1}{\e_j}\Bigl(v_j^{(k)}\Bigl(\frac{t}{\e_j}\Bigr)\Bigr)^2\Big)J(y,t)\,dt\,d\mathcal{H}^{d-1}(y) \\
    &=& \int_{M\cap\Omega}\int_{-\frac1{\sqrt{\e_j}}}^{\frac1{\sqrt{\e_j}}}\bigl(W(v_j(s))+(v_j^{(k)}(s))^2\bigr)J(y,s\e_j)\,ds\,d\mathcal{H}^{d-1}(y) \\
    &\leq & \bigg(\sup_{\substack{{y\in M}\\t\in(-\sqrt{\e_j},\sqrt{\e_j})}}J(y,t)\bigg)\int_{M\cap\Omega}\int_{-\infty}^{+\infty}(W(v_j(s))+(v_j^{(k)}(s))^2)ds\,d\mathcal{H}^{d-1}(y).
\end{eqnarray*}
Since $M$ is compact, $J(y,t)$ converges uniformly to $1$ as $t\to0$; hence, by \eqref{eq:vj is an optimal profile limd} we obtain
\[
\limsup_{j\to+\infty}  I_j^k \leq m_k \mathcal{H}^{d-1}(M\cap\Omega) = m_k  \textup{Per}(E;\Omega) = m_k \mathcal{H}^{d-1}(S(u)).
\]

We now prove that for fixed $\ell\in\{1,...,k-1\}$ we have 
$\limsup\limits_{j\to+\infty}I_j^{\ell}=0$. 
To this end, we set $w_j(t):=v_j(t/\e_j)$ and $A_j:=(-\sqrt{\e_j},\sqrt{\e_j})$, and write
\begin{equation}\label{limsupclaim1}
    I_j^\ell=\e_j^{2k-1-2\ell}\int_{A_j}\Bigl(v_j^{(\ell)}\Bigl(\frac{t}{\e_j}\Bigr)\Bigr)^2\,dt = \e_j^{2k-1} \int_{A_j}(w_j^{(\ell)}(t))^2\,dt.
\end{equation}
By a change of variable, we have
\begin{equation*}
    \int_{A_j}\Big(\frac{1}{\e_j}W(w_j)+\e_j^{2k-1}(w_j^{(k)})^2\Big)dt=\int_{-\infty}^{+\infty}\bigl(W(v_j)+(v_j^{(k)})^2\bigr)\,dt;
\end{equation*}
hence, by \eqref{eq:vj is an optimal profile limd} and by the fact that $W\geq0$, we infer that
\begin{equation}\label{boundw}
    \|w_j^{(k)}\|_{L^2(A_j)}\leq C \e_j^{\frac{1-2k}{2}}.
\end{equation}
Note that for every positive integer $j$ and for any $h\in\{0,...,k-1\}$ there exists $z^h_j\in [-\sqrt{\e_j},\sqrt{\e_j}]$ such that $w^{(h)}_j(z_j^h)=0$. Indeed, this holds with $z^h_j=\sqrt{\e_j}$ for $h\in\{1,...,k-1\}$, while the existence of $z^0_j$ follows by the continuity of $w_j$. Hence, applying the Fundamental Theorem of Calculus and using \eqref{boundw}, we get
\begin{eqnarray*}
    |w_j(t)|\leq \e_j^{\frac{2k-1}{4}}\|w_j^{(k)}\|_{L^2(A_j)}\leq C \e_j^{\frac{1-2k}{4}} \quad \text{for every } t\in A_j,
\end{eqnarray*}
so that 
\begin{equation}\label{boundw2}
    \|w_j\|_{L^2(A_j)}\leq C \e_j^{\frac{1-k}{2}}.
\end{equation}
By the interpolation inequality \eqref{interpolation1}, it holds
\begin{equation*}
\|w_j^{(\ell)}\|_{L^2(A_j)}\leq C\Bigl(\|w_j^{(k)}\|_{L^2(A_j)}^{\frac{\ell}{k}}\|w_j\|_{L^2(A_j)}^{\frac{k-\ell}{k}}+|A_j|^{-\ell}\|w_j\|_{L^2(A_j)}\Bigr). 
\end{equation*}
Hence, combining \eqref{boundw} and \eqref{boundw2} we get
\begin{equation}\label{eq:normal2derivataw}
    \|w_j^{(\ell)}\|_{L^2(A_j)}\leq C \Bigl(\e_j^{(\frac{1-2k}{2})(\frac{\ell}{k})} \e_j^{(\frac{1-k}{2})(\frac{k-\ell}{k})}+\e_j^{-\frac{\ell}{2}+\frac{1-k}{2}}\Bigr)  = C \e_j^{\frac{-k-\ell+1}{2}}.
\end{equation}
Note that the above formula, up to a constant factor that is independent of $j$, reduces to \eqref{boundw} when $\ell=k$.
Recalling \eqref{limsupclaim1} we conclude 
\begin{equation*}
    I_j^\ell = \e_j^{2k-1}\|w_j^{(\ell)}\|_{L^2(A_j)}^2 \leq C \e^{k-\ell},
\end{equation*}
which tends to 0 since $\ell<k$.

Finally, we prove that for fixed $r<\ell$, we have
$\limsup\limits_{j\to+\infty}I^{r,\ell}_j=0$.

Recalling that $w_j(t)=v_j(t/\e_j)$, by H{\"o}lder's inequality we have
\begin{eqnarray*}
    I_{j}^{r,\ell}&=& \e_j^{2k-1-r-\ell}\int_{A_j}{}\Bigr|v_j^{(r)}\Bigl(\frac{t}{\e_j}\Bigr)\Bigr|\Bigl|v_j^{(\ell)}\Bigl(\frac{t}{\e_j}\Bigr)\Bigr|\,dt\\
     &=& \e_j^{2k-1}\int_{A_j}|w_j^{(r)}(t)||w_j^{(\ell)}(t)|\,dt\\
     &\leq&\e_j^{2k-1}\|w_j^{(r)}\|_{L^2(A_j)}\,\|w_j^{(\ell)}\|_{L^2(A_j)}.
\end{eqnarray*}
Since \eqref{eq:normal2derivataw} holds for every $\ell\in\{1,...,k\}$, we get
\begin{equation*}
    I_j^{r,\ell} \leq C \e_j^{2k-1} \e_j^{\frac{-k-r+1}{2}} \e_j^{\frac{-k-\ell+1}{2}} = C \e_j^{k-\frac{r+\ell}{2}},
\end{equation*}
which tends to $0$ since $r+\ell < 2k$, concluding the proof.
\end{proof}

\noindent{\bf Acknowledgements.}
 This paper is based on work supported by the National Research Project PRIN 2022J4FYNJ  ``Variational methods for stationary and evolution problems with singularities and interfaces'' 
 funded by the Italian Ministry of University and Research and by the GNAMPA Project ``Asymptotic analysis of nonlocal variational problems'' funded by INdAM.  
The authors are members of GNAMPA of INdAM. 

\bibliographystyle{abbrv}
\bibliography{References}

\begin{thebibliography}{10}

\bibitem{alberti1998phase}
G.~Alberti, G.~Bouchitt{\'e}, and P.~Seppecher.
\newblock Phase transition with the line-tension effect.
\newblock {\em Arch. Rational Mech. Analysis}, 144:1--46, 1998.

\bibitem{ABG}
R.~Alicandro, A.~Braides, and M.~S. Gelli.
\newblock Free-discontinuity problems generated by singular perturbation.
\newblock {\em Proc. Royal Soc. Edinburgh}, 128A:1115--1129, 1998.

\bibitem{Bach}
A.~Bach.
\newblock Anisotropic free-discontinuity functionals as the {$\Gamma$}-limit of
  second-order elliptic functionals.
\newblock {\em ESAIM Control Optim. Calc. Var.}, 24(3):1107--1140, 2018.

\bibitem{bcg}
G.~Bellettini, A.~Chambolle, and M.~Goldman.
\newblock The {$\Gamma$}-limit for singularly perturbed functionals of
  {P}erona-{M}alik type in arbitrary dimension.
\newblock {\em Mathematical Models and Methods in Applied Sciences},
  24(06):1091--1113, 2014.

\bibitem{BBL}
M.~Bonnivard, E.~Bretin, and A.~Lemenant.
\newblock Numerical approximation of the {S}teiner problem in dimension 2 and
  3.
\newblock {\em Math. Comp.}, 89(321):1--43, 2020.

\bibitem{BDS}
G.~Bouchitt{\'e}, C.~Dubs, and P.~Seppecher.
\newblock Regular approximation of free-discontinuity problems.
\newblock {\em Mathematical Models and Methods in Applied Sciences},
  10(7):1073--1097, 2000.

\bibitem{BLN98}
A.~Braides.
\newblock {\em Approximation of Free-Discontinuity Problems}.
\newblock Springer-Verlag, Berlin, 1998.

\bibitem{BMS}
A.~Braides, M.~Maslennikov, and L.~Sigalotti.
\newblock Homogenization by blow-up.
\newblock {\em Applicable Anal.}, 87:1341--1356, 2008.

\bibitem{ChCo}
M.~Chermisi and S.~Conti.
\newblock Multiwell rigidity in nonlinear elasticity.
\newblock {\em SIAM J. Math. Anal.}, 42(5):1986--2012, 2010.

\bibitem{CDMFL}
M.~Chermisi, G.~Dal~Maso, I.~Fonseca, and G.~Leoni.
\newblock Singular perturbation models in phase transitions for second-order
  materials.
\newblock {\em Indiana Univ. Math. J.}, 60(2):367--409, 2011.

\bibitem{CFL}
S.~Conti, I.~Fonseca, and G.~Leoni.
\newblock A {$\Gamma$}-convergence result for the two-gradient theory of phase
  transitions.
\newblock {\em Comm. Pure Appl. Math.}, 55(7):857--936, 2002.

\bibitem{Cos}
S.~Conti and B.~Schweizer.
\newblock A sharp-interface limit for a two-well problem in geometrically
  linear elasticity.
\newblock {\em Arch. Ration. Mech. Anal.}, 179(3):413--452, 2006.

\bibitem{FLL}
I.~Fonseca, P.~Liu, and X.~Y. Lu.
\newblock Higher order {A}mbrosio-{T}ortorelli scheme with non-negative
  spatially dependent parameters.
\newblock {\em Adv. Calc. Var.}, 16(4):885--902, 2023.

\bibitem{Fonseca2000}
I.~Fonseca and C.~Mantegazza.
\newblock Second order singular perturbation models for phase transitions.
\newblock {\em SIAM J. Math. Anal.}, 31:1121--1143, 2000.

\bibitem{FM}
I.~Fonseca and S.~M\"uller.
\newblock Quasiconvex integrands and lower semicontinuity in ${L}^1$.
\newblock {\em SIAM J. Math. Anal.}, 23:1081--1098, 1992.

\bibitem{gp}
M.~Gobbino and N.~Picenni.
\newblock A quantitative variational analysis of the staircasing phenomenon for
  a second order regularization of the {P}erona-{M}alik functional.
\newblock {\em Transactions of the American Mathematical Society}, 2023.

\bibitem{Gurtin}
M.~E. Gurtin.
\newblock Some results and conjectures in the gradient theory of phase
  transitions.
\newblock In {\em Metastability and {I}ncompletely {P}osed {P}roblems}, pages
  135--146. Springer, 1987.

\bibitem{Ho}
P.~Hornung.
\newblock A {$\Gamma$}-convergence result for thin martensitic films in
  linearized elasticity.
\newblock {\em SIAM J. Math. Anal.}, 40(1):186--214, 2008.

\bibitem{IZ}
R.~Ignat and H.~Zorgati.
\newblock Dimension reduction and optimality of the uniform state in a
  phase-field-crystal model involving a higher-order functional.
\newblock {\em J. Nonlinear Sci.}, 30(1):261--282, 2020.

\bibitem{leoni}
G.~Leoni.
\newblock {\em A First Course in Sobolev Spaces}.
\newblock American Mathematical Society, Providence, 2017.

\bibitem{Maggi}
F.~Maggi.
\newblock {\em Sets of Finite Perimeter and Geometric Variational Problems: an
  Introduction to Geometric Measure Theory}.
\newblock Cambridge University Press, Cambridge, 2012.

\bibitem{MM}
L.~Modica and S.~Mortola.
\newblock Un esempio di {$\Gamma$}-convergenza.
\newblock {\em Boll. Un. Mat. It. B}, 14:285--299, 1977.

\bibitem{mo}
M.~Morini.
\newblock Sequences of singularly perturbed functionals generating
  free-discontinuity problems.
\newblock {\em SIAM Journal on Mathematical Analysis}, 35(3):759--805, 2003.

\bibitem{Negri}
M.~Negri.
\newblock {$\Gamma$}-convergence for high order phase field fracture: continuum
  and isogeometric formulations.
\newblock {\em Comput. Methods Appl. Mech. Engrg.}, 362:112858, 21, 2020.

\bibitem{Solci2024}
M.~Solci.
\newblock Free-discontinuity problems generated by higher-order singular
  perturbations, 2024, https://arxiv.org/abs/2402.10656.

\end{thebibliography}

\end{document}